\documentclass[11pt]{amsart}
\usepackage[colorlinks=true, pdfstartview=FitV, linkcolor=blue, citecolor=blue, urlcolor=blue, breaklinks=true]{hyperref}
\usepackage{amsmath,amsfonts,amssymb,amsthm,amscd,comment,euscript}
\usepackage{xcolor}
\usepackage{paralist}
\usepackage[vcentermath]{youngtab}
\usepackage[arrow, matrix, curve]{xy}
\usepackage[latin1]{inputenc}
\usepackage{tikz}
\usepackage{array}
\newtheorem{cor}{Corollary}[subsection]
\newtheorem{lem}{Lemma}[subsection]
\newtheorem{prop}{Proposition}[subsection]

\newcommand{\nc}{\newcommand}

\newcounter{cnt}
 \makeatletter
\def\mydggeometry{\makeatletter\dg@YGRID=1\dg@XGRID=20\unitlength=0.003pt\makeatother}\makeatother \theoremstyle{remark}
\numberwithin{equation}{section}

\theoremstyle{definition}
\newtheorem*{defn}{Definition}
\newtheorem*{ex}{Example}
\theoremstyle{definition}
\newtheorem{thm}{Theorem}[subsection]
\newtheorem*{theom}{Theorem}

\newtheorem{rem}{Remark}[subsection]

\newcommand\Lg{\mathfrak{g}}
\newcommand\Lh{\mathfrak{h}}
\newcommand\Lb{\mathfrak{b}}
\newcommand\Ln{\mathfrak{n}}
\newcommand\U{\mathfrak{U}}

\nc{\C}{\mathbb C }
\nc{\D}{\mathbb D }
\nc{\Z}{\mathbb Z }
\nc{\N}{\mathbb N }
\nc{\R}{\mathbb R }
\nc{\Q}{\mathbb Q }
\newcommand{\wt}{\operatorname{wt}}

\begin{document}

\title[Crystal bases as tuples of integer sequences]{Crystal bases as tuples of integer sequences}
\author{Deniz Kus}
\address{Deniz Kus:\newline
Mathematisches Institut, Universit\"at zu K\"oln, Germany}
\email{dkus@math.uni-koeln.de}
\thanks{The author was sponsored by the "SFB/TR 12 - Symmetries and Universality in Mesoscopic Systems".}

\begin{abstract}
We describe a set $\mathcal{R}^{\infty}$ consisting of tuples of integer sequences and provide certain explicit maps on it. We show that this defines a semiregular crystal for $\mathfrak{sl}_{n+1}$ and $\mathfrak{sp}_{2n}$ respectively. Furthermore we define for any dominant integral weight $\lambda$ a connected subcrystal $\mathcal{R}(\lambda)$ in $\mathcal{R}^{\infty}$, such that this crystal is isomorphic to the crystal graph $B(\lambda)$. Finally we provide an explicit description of these connected crystals $\mathcal{R}(\lambda)$.
\end{abstract}
\maketitle \thispagestyle{empty}
\section{Introduction}
Let $\Lg$ be a symmetrizable Kac-Moody algebra and $\U_q(\Lg)$ be the corresponding quantum algebra. For these quantum algebras Kashiwara developed the crystal bases theory for integrable modules in \cite{Kashi91} and thus provided a remarkable combinatorial tool for studying these modules. In particular crystal bases can be viewed as bases at $q=0$ and they contain structures of edge-colored oriented graphs satisfying a set of axioms, called the crystal graphs. These crystal graphs have certain nice properties, for instance characters of $\U_q(\Lg)$-modules can be computed and the decomposition of tensor
products of modules into irreducible ones can also be determined from the crystal graphs, to name just a few. It is thus an important problem to have explicit realizations of crystal graphs.\par
There are many such realizations, combinatorial and geometrical, worked out during the last decades, for instance we refer to (\cite{NK94},\cite{L94},\cite{NS97},\cite{Nak01}). In \cite{NK94} the authors give a tableaux realization of crystal graphs for irreducible modules over the quantum algebra for all classical Lie algebras, which is a purely combinatorial model. Another significant combinatorial model for any symmetrizable Kac-Moody algebra is provided in \cite{L94}, called Littelmann's path model. The underlying set here is a set of piecewise linear maps, and the crystal graph of an irreducible module of any dominant integral highest weight $\lambda$ can be generated by an algorithm using the straight path
connecting 0 and $\lambda$.\par
A geometrical realization of crystals is also known and is provided by Nakajima \cite{Nak01} by showing that there exists a crystal structure on the set of irreducible components of a lagrangian subvariety of the quiver variety $\mathcal{M}$. This realization can be translated into a purely combinatorial model, the set of Nakajima monomials, where the action of the Kashiwara operators can be understood as a multiplication with monomials. Moreover, it is shown in \cite{Kashi03} that the connected component of any highest weight monomial of highest weight $\lambda$ is isomorphic to the crystal graph $B(\lambda)$ obtained from Kashiwara's crystal bases theory. For special highest weight monomials these connected components are explicitly characterized for $\mathfrak{sl}_{n+1}$ in \cite{KKS04} and for the other classical Lie algebras in \cite{KKS2004}. A combinatorial isomorphism from connected components corresponding to arbitrary highest weight monomials of highest weight $\lambda$ and those in \cite{KKS04},\cite{KKS2004} is provided in \cite{M12} for the types $A$ and $C$.\par
In this paper we introduce a set $\mathcal{R}^{\infty}$ consisting of tuples of integer sequences, i.e. a typical element in $\mathcal{R}^{\infty}$ is given by $$\mathbf x=(x_1,x_2,\cdots)\in \mathcal{R}^{\infty},$$ where each component $x_j$ consists of certain ordered pairs of integers, $x_j=(i_1,i'_1)\cdots(i_s,i'_s)$ (see Definition~\ref{maindef}). Furthermore, the number of non-zero components is finite.
We provide certain maps on $\mathcal{R}^{\infty}$, the Kashiwara operators $\tilde{e}_l$, $\tilde{f}_l$ and maps $\epsilon_l$, $\varphi_l$ for all $l=1,\ldots,n$ and prove that $\mathcal{R}^{\infty}$ is a semiregular crystal if $\Lg$ is $\mathfrak{sl}_{n+1}$ or $\mathfrak{sp}_{2n}$ (see Definition~\ref{abcrystal} and Proposition~\ref{crystalprop}).\par
Moreover, we introduce for any dominant integral weight $\lambda$ a subcrystal $\mathcal{R}(\lambda)$ as the connected component of $\mathcal{R}^{\infty}$ containing a highest weight element $r_{\lambda}$ and prove the following theorem:
\begin{theom}
Let $\lambda$ be a dominant integral weight, then there exists a crystal isomorphism $$\mathcal{R}(\lambda)\longrightarrow B(\lambda),$$ mapping $r_{\lambda}$ to the highest weight element $b_\lambda\in B(\lambda).$
\end{theom}
Therefore, similar to the setting of Nakajima monomials, a natural question arises, namely can one characterize for each dominant integral weight $\lambda$ explicitly the sequences appearing in $\mathcal{R}(\lambda)$?
We answer this question by describing explicitly these connected components (for the special linear Lie algebra in Theorem~\ref{char} and the symplectic Lie algebra in Theorem~\ref{char2}).
\vskip 5pt 
Our paper is organized as follows: in Section~\ref{section2} we fix some notation and review briefly the crystal theory. In Section~\ref{section3} we present the main definitions, especially the definition of $\mathcal{R}^{\infty}$ and we equip our main object with a crystal structure. In Section~\ref{section4} Nakajima monomials are recalled. In Section~\ref{section5} we introduce for any dominant integral weight $\lambda$ the subcrystals $\mathcal{R}(\lambda)$ and describe them explicitly. Finally, in Section~\ref{section6} we prove that they are isomorphic to $B(\lambda)$.
\vskip12pt
%
\section{Notations and a review of crystal theory}\label{section2}
Let $\Lg$ be a complex simple Lie algebra of rank $n$ with index set $I=\{1,\cdots,n\}$ and fix a Cartan subalgebra $\Lh$ in $\Lg$ and a Borelsubalgebra $\Lb\supseteq\Lh$. We denote by $\Phi\subseteq \Lh^*$ the root system of the Lie algebra, and, corresponding to the choice of $\Lb$ let $\Phi^+$ be the subset of positive roots. Further, we denote by $\Pi=\{\alpha_1,\cdots,\alpha_n\}$ the corresponding basis of $\Phi$ and the basis of the dual root system $\Phi^{\vee}\subseteq \Lh$ is denoted by $\Pi^\vee=\{\alpha_1^\vee,\cdots,\alpha_n^\vee\}$. Let  $\Lg=\Ln^+\oplus\Lh\oplus \Ln^-$ be a Cartan decomposition and for a given root $\alpha\in \Phi$ let $\Lg_{\alpha}$ be the corresponding root space.
For a dominant integral weight $\lambda$ we denote by $V(\lambda)$ the irreducible $\Lg$-module with
highest weight $\lambda$. Fix a highest weight vector $v_\lambda\in V(\lambda)$,
then $V(\lambda)=\U(\Ln^-)v_\lambda$, where $\U(\Ln^-)$ denotes the universal
enveloping algebra of $\Ln^-$. For an indetermined element $q$ we denote by $\U_q(\Lg)$ be the corresponding quantum algebra. The theory of studying modules of quantum algebras is quite parallel to that of Kac-Moody algebras and the irreducible modules are classified again in terms of highest weights (see \cite{HK02}). Using the crystal bases theory, introduced by Kashiwara in \cite{Kashi91}, we can compute the character of an integrable module $M$ in the category $\mathcal{O}^q$ as follows:
$$ ch M=\sum_{\mu}\sharp (B_{\mu}) e^{\mu},$$ whereby $(L,B)$ is the crystal bases of $M$ (see \cite{HK02}). The crystal graph associated to the irreducible module of highest weight $\lambda$ is denoted by $B(\lambda)$. So finding expressions for the characters can be achieved by finding explicit combinatorial description of crystal bases. For some examples we refer to (\cite{NK94},\cite{L94},\cite{NS97}).\par
From now on we assume that $\Lg$ is a classical Lie algebra of type $A_n$ or $C_n$. Note that the positive roots are all of the following form $$\mbox{Type } A_n: \alpha_{i,j}=\alpha_i+\alpha_{i+1}+\cdots+\alpha_j,\ \mbox{for $1\leq i\leq j\leq n$}$$ 
$$ \mbox{Type } C_n:  \alpha_{i,j}=\alpha_i+\alpha_{i+1}+\cdots+\alpha_j,\ \mbox{for $1\leq i\leq j\leq n$}$$$$ \alpha_{i,\overline{j}}=\alpha_i+\alpha_{i+1}+\cdots+\alpha_n+\alpha_{n-1}+\cdots+\alpha_j,\ \mbox{for $1\leq i\leq j\leq n$}.$$
Furthermore let $P=\bigoplus_{i\in I}\Z\omega_i$ be the set of classical integral weights and $P^+=\bigoplus_{i\in I}\Z_{+}\omega_i$ be the set of classical dominant integral weights. Before we discuss the crystal bases theory in detail we review first the notion of abstract crystals.
%
\subsection{Abstract crystals}
Crystal bases of integrable $\U_q(\Lg)$-modules in the category $\mathcal{O}^q$ are characterized by certain maps satisfying some properties. One can define the abstract notion of crystals associated with a Cartan datum as follows:
\begin{defn}\label{abcrystal}
Let $I$ be a finite index set and let $A=(a_{i,j})_{i,j\in I}$ be a generalized Cartan matrix with the Cartan datum $(A,\Pi,\Pi^{\vee},P,P^{\vee})$. A crystal associated with the Cartan datum $(A,\Pi,\Pi^{\vee},P,P^{\vee})$ is a set $B$ together with maps $\wt:B \rightarrow P$, $\tilde{e}_l,\tilde{f}_l: B\rightarrow B \cup \{0\}$, and $\epsilon_l,\varphi_l:B\rightarrow \Z\cup \{-\infty\}$ satisfying the following properties for all $l\in I$:
\begin{enumerate}
\item $\varphi_l(b)=\epsilon_l(b)+\langle \alpha_l^\vee,\wt(b)\rangle$
\item $\wt(\tilde{e}_{l}b)=\wt(b)+\alpha_l$ if $\tilde{e}_{l}b\in B$
\item $\wt(\tilde{f}_{l}b)=\wt(b)-\alpha_l$ if $\tilde{f}_{l}b\in B$
\item $\epsilon_l(\tilde{e}_{l}b)=\epsilon_l(b)-1$, $\varphi_l(\tilde{e}_{l}b)=\varphi_l(b)+1$ if $\tilde{e}_{l}b\in B$
\item $\epsilon_l(\tilde{f}_{l}b)=\epsilon_l(b)+1$, $\varphi_l(\tilde{f}_{l}b)=\varphi_l(b)-1$ if $\tilde{f}_{l}b\in B$
\item $\tilde{f}_{l}b=b'$ if and only if $\tilde{e}_{l}b'=b$ for $b,b'\in B$
\item if $\varphi_l(b)=-\infty$ for $b\in B$, then $\tilde{f}_{l}b=\tilde{e}_{l}b$=0.
\end{enumerate}
Furthermore a crystal $B$ is said to be semiregular if the equalities $$\epsilon_l(b)=\max\{k\geq0|\tilde{e}_{l}^{k}b\neq 0\},\quad \varphi_l(b)=\max\{k\geq0|\tilde{f}_{l}^{k}b\neq 0\}$$ hold.
\end{defn}
The maps $\tilde{e}_l$ and $\tilde{f}_l$ are called Kashiwara's crystal operators and the map $\wt$ is called the weight function. 
So, on the one hand one can associate to any integrable $\U_q(\Lg)$-module a set $B$ satisfying the properties from Definition~\ref{abcrystal} and on the other hand one can study the notion of abstract crystals. A natural question which arises at this point is therefore the following: can one determine whether an abstract crystal is the crystal of a module? Stembridge \cite{S2003} gave a set of local axioms to characterize the set of crystals of module in the class of all crystals when $\Lg$ is simply-laced and a list of local axioms for $B_2$-crystals is provided in \cite{DKK07}. In the following sections we define our underlying set and realize the crystal obtained from Kashiwara's crystal bases theory for the types $A_n$ and $C_n$. We start by equipping our underlying set with an abstract crystal structure and later we prove that this crystal is the crystal of a module.
%
\section{Tuples of integer sequences as crystals}\label{section3}
In this section we introduce a set $\mathcal{R}^{\infty}$ consisting of tuples of integer sequences (see Definition~\ref{maindef}) and a crystal structure on it in the sense of Definition~\ref{abcrystal}. Our purpose is to identify for any dominant integral weight $\lambda$ certain subcrystals $\mathcal{R}(\lambda)$, i.e. $\bigcup_{\lambda\in P^+}\mathcal{R}(\lambda)\subseteq \mathcal{R}^{\infty}$, such that $\mathcal{R}(\lambda)$ has a strong connection to the crystal graph $B(\lambda)$ (see Corollary~\ref{mainthm1}). 
%
\subsection{Set of tuples of integer sequences}
In order to define $\mathcal{R}^{\infty}$ we consider a total order on $\mathbf I=\{1,\cdots,n\}$  if $\Lg$ is of type $A_n$ and a total order on $\mathbf I=\{1,\cdots,n,\overline{n-1},\cdots,\overline{1}\}$ if $\Lg$ is of type $C_n$, namely $$1< 2 < \cdots < n$$ and $$1< 2 < \cdots < n< \overline{n-1}<\cdots<\overline{1},$$ respectively. Furthermore, especially in Section~\ref{section5}, we need for type $C_n$ the following bijective map $$\bar{}\ :\mathbf I\longrightarrow \mathbf I$$ $$n\mapsto n,\ \overline{i}\mapsto i,\ i\mapsto\overline{i},\mbox{ for } i\in\{1,\cdots,n-1\}.$$\par
For $1\leq i\leq n$,$\ s\in\Z_{\geq 0}$ we set $\mathcal{R}^i_s$ to be the set of all sequences $(i_1,i'_1)\cdots(i_s,i'_s)$ with $i_j,i'_j\in\mathbf I$, such that \begin{eqnarray}\label{cond}&1\leq i_s< i_{s-1}< \cdots < i_1 \leq i \leq i'_1< i'_{2}\cdots < i'_s \leq \max\mathbf I
&\\& \notag i'_j\leq \overline{i_j},\ j=1,\ldots,s,\end{eqnarray} where $\max\mathbf I$ denotes the maximal element in $\mathbf I$ with respect to $<$. We denote by $\emptyset_i$ the unique element in $\mathcal{R}^i_0$. 
\begin{defn}\label{maindef}
We define $\mathcal{R}^{\infty}$ to be the set of all infinite sequences $\mathbf x=(x_1 ,x_2 ,x_3,\cdots)$ where each component $x_j$ is contained in $\mathcal{R}=\bigcup_{s,i}\mathcal{R}^i_s \cup \{0\}$ and only finitely many components are non-zero. We identify $\mathcal{R}^k$ with the sequences of the form $(x_1 ,\cdots,x_k,0,0,\cdots)$.
\end{defn}
Before we mention the crystal structure on $\mathcal{R}^{\infty}$ we will initially introduce a list of properties. We need these to define the Kashiwara operators. Let $x=(i_1,i'_1)\cdots(i_s,i'_s)\in\mathcal{R}^i_s$ be an arbitrary element and fix $l\in I$:\vskip 5pt
\label{conditions}
\begin{enumerate}[(a)]
\item \begin{align*}l\notin\{i_1,\cdots,i_{s},i'_1,\cdots,i'_{s}\} &\mbox{ and } 
\begin{cases}\mbox{$l+1\in\{i_1,\cdots,i_{s}\}$},& \text{if $l<i$}\\
\mbox{$l-1\in\{i'_1,\cdots,i'_{s}\}$},& \text{if $l>i$}\\
\mbox{/},& \text{if $l=i$}\end{cases}
&\\& \mbox{ and } (\overline{l+1}\in\{i'_1,\cdots,i'_{s}\} \vee\ \overline{l}\notin\{i'_1,\cdots,i'_{s}\})\end{align*}
\end{enumerate}
\vskip12pt
\begin{enumerate}[(a')]
\item Replace in (a) $\vee$ by $\wedge$
\end{enumerate}
\vskip12pt
\begin{enumerate}[(b)]
\item  \begin{align*}\overline{l}\notin\{i'_1,\cdots,i'_{s}\}, \overline{l+1}\in\{i'_1,&\cdots,i'_{s}\} &\\&\mbox{ and }\begin{cases}\mbox{$l+1\notin\{i_1,\cdots,i_{s}\}\vee l\in\{i_1,\cdots,i_{s}\}$},& \text{if $l<i$}\\
\mbox{$l-1\notin\{i'_1,\cdots,i'_{s}\}\vee l\in\{i'_1,\cdots,i'_{s}\}$},& \text{if $l>i$}\\
\mbox{$l=i_1 \vee l=i'_1$},& \text{if $l=i$}\end{cases}\end{align*}
\end{enumerate}
\vskip12pt
\begin{enumerate}[(c)]
\item \begin{align*}l\in\{i_1,\cdots,i_{s},i'_1,\cdots,i'_{s}\} &\mbox{ and } 
\begin{cases}\mbox{$l+1\notin\{i_1,\cdots,i_{s}\}$},& \text{if $l<i$}\\
\mbox{$l-1\notin\{i'_1,\cdots,i'_{s}\}$},& \text{if $l>i$}\\
\mbox{$(i_1,i'_1)=(l,l)$},& \text{if $l=i$}\end{cases}
&\\& \mbox{ and } (\overline{l+1}\in\{i'_1,\cdots,i'_{s}\} \vee\ \overline{l}\notin\{i'_1,\cdots,i'_{s}\})\end{align*}
\end{enumerate}
\vskip12pt
\begin{enumerate} [(d)]
\item \begin{align*}\overline{l}\in\{i'_1,\cdots,i'_{s}\}, \overline{l+1}\notin\{i'_1,&\cdots,i'_{s}\} &\\&\mbox{ and }\begin{cases}\mbox{$l+1\notin\{i_1,\cdots,i_{s}\}\vee l\in\{i_1,\cdots,i_{s}\}$},& \text{if $l<i$}\\
\mbox{$l-1\notin\{i'_1,\cdots,i'_{s}\}\vee l\in\{i'_1,\cdots,i'_{s}\}$},& \text{if $l>i$}\\
\mbox{$l=i_1 \vee l=i'_1$},& \text{if $l=i$}\end{cases}\end{align*}
\end{enumerate}
\vskip12pt
\begin{enumerate}[(d')]
\item Replace everywhere in (d) $\vee$ by $\wedge$.
\end{enumerate}
Let us consider an example.
\begin{ex}\mbox{}
\begin{enumerate}
\item Let $\Lg=A_5$, $l=2$ and $$x_1=(3,4)(1,5)\in\mathcal{R}_2^3,\ x_2=(2,3)\in\mathcal{R}_1^2,$$ then $x_1$ satisfies (a) while $x_2$ violates (a).
\item Let $\Lg=C_3$, $l=1$ and $$x=(1,\overline{2})\in\mathcal{R}_1^2,$$ then $x$ satisfies (c).
\end{enumerate}
\end{ex}
\begin{rem}
If $x$ satisfies (a') and (d') respectively, then it satisfies also (a) and (d) respectively. If $\Lg$ is further of type $A_n$, these properties can be simplified. In particular, the properties (a'), (b), (d) and (d') are superfluous.
\end{rem}
Henceforth we define a crystal structure on $\mathcal{R}^{\infty}$, such that the semiregularity holds. For this let $\mathbf x=(x_1 ,x_2 ,x_3,\cdots)$ be such a sequence with finitely many components different from zero; recall that each component is a sequence as in (\ref{cond}). The weight function is given by \begin{equation}\label{weight}\wt(\mathbf x)=\sum^n_{i=1}c_i\omega_i-\sum_j \wt(x_j),\end{equation} 
where $$\wt(x_j)=\begin{cases}\sum^s_{j=1} \alpha_{i_j,i'_j} &\text{if $x_j=(i_1,i'_1)\cdots(i_s,i'_s)$}\\
0& \text{if $x_j=0$ or $x_j\in\{\emptyset_1,\cdots,\emptyset_n\}$} \end{cases}$$ and $c_i=\sharp\{x_j\neq 0|x_j\in\bigcup_s\mathcal{R}_s^i\}.$ Suppose that the non-zero components in $\mathbf x$ are given by $x_{q_1}\in\mathcal{R}_{s_1}^{j_1},\cdots,x_{q_k}\in \mathcal{R}_{s_k}^{j_k}$. 
For fixed $l\in I$ we define the following maps:\vskip3pt

$\bullet$ for $2\leq j\leq k+1$ let $\sigma^j_l:\mathcal{R}^{\infty}\rightarrow \Z_{\geq0}$ be the map given by $$\sigma^j_l(\mathbf x)=a^j_l(\mathbf x)+b^j_l(\mathbf x),$$ where $$a^j_l(\mathbf x)=\sharp\{x_{q_p}|1\leq p\leq j-1, x_{q_p} \mbox{ satisfies (a) or (b)}\}$$
$$b^j_l(\mathbf x)=\sharp\{x_{q_p}|1\leq p\leq j-1, x_{q_p}\mbox{ satisfies (a')}\}.$$ 
\vskip5pt
Furthermore we define $\theta_l(x_{q_p})$ to be the sequence which arises out of $x_{q_p}$ by $$\begin{cases}\mbox{replacing $l+1$  by $l$},&\text{if $l<j_p$}\\
\mbox{replacing $l-1$ by $l$},&\text{if $l>j_p$}\\
\mbox{adding $(l,l)$},&\text{if $l=j_p$,}\end{cases}$$
if $x_{q_p}$ satisfies (a). If $x_{q_p}$ satisfies (b) let $\theta_l(x_{q_p})$ be the sequence which arises out of $x_{q_p}$ by replacing $\overline{l+1}$ by $\overline{l}$. If neither (a) nor (b) is fulfilled, we set $\theta_l(x_{q_p})=0$.

$\bullet$ For $2\leq j\leq k$ let $\tau^j_l:\mathcal{R}^{\infty}\rightarrow \Z_{\geq0}$ be the map given by $$\tau^j_l(\mathbf x)=c^j_l(\mathbf x)+d^j_l(\mathbf x),$$ where 
$$c^j_l(\mathbf x)=\sharp\{x_{q_p}|2\leq p\leq j, x_{q_p} \mbox{ satisfies (c) or (d)}\}$$
$$d^j_l(\mathbf x)=\sharp\{x_{q_p}|2\leq p\leq j, x_{q_p}\mbox{ satisfies (d')}\}.$$ 
\vskip5pt
We define $\rho_l(x_{q_p})$ to be the sequence which arises out of $x_{q_p}$ by $$\begin{cases}\mbox{replacing $l$ by $l+1$},& \text{if $l<j_p$}\\
\mbox{replacing $l$ by $l-1$},&\text{if $l>j_p$}\\
\mbox{erasing $(l,l)$},&\text{if $l=j_p$},\end{cases}
$$ 
if $x_{q_p}$ satisfies (c). If $x_{q_p}$ satisfies (d) let $\rho_l(x_{q_p})$ be the sequence which arises out of $x_{q_p}$ by replacing $\overline{l}$ by $\overline{l+1}$. If neither (c) nor (d) is fulfilled, we set $\rho_l(x_{q_p})=0$.

\begin{rem}\label{remark1}\mbox{}
\begin{enumerate}
\item For $j=1$ we set $\sigma^j_l=\tau^j_l=0$.
\item Note that the image of $x\in\bigcup_s\mathcal{R}^i_s$ under the map $\theta_l$ and $\rho_l$ respectively is contained in $\bigcup_s\mathcal{R}^i_s\cup \{0\}$, i.e. $$\theta_l:\bigcup_s\mathcal{R}^i_s\longrightarrow\bigcup_s\mathcal{R}^i_s\cup \{0\}$$
$$\rho_l:\bigcup_s\mathcal{R}^i_s\longrightarrow\bigcup_s\mathcal{R}^i_s\cup \{0\}$$
\end{enumerate}
\end{rem}
One important fact about these maps is described in the next lemma.
\begin{lem}\label{vor1}
Let $x,\tilde{x}$ be non-zero sequences as in (\ref{cond}), then we have $$\theta_l(x)=\tilde{x}\mbox{ if and only if } \rho_l(\tilde{x})=x.$$
\proof
One can easily show that $x$ satisfies (a) if and only if $\tilde{x}$ satisfies (c). Hence we can suppose that $x$ does not fulfill all properties enumerated in (a). By observing the action we see that $\tilde{x}$ arises from $x$ by replacing $\overline{l+1}$ by $\overline{l}$, which means that (c) is violated. In particular $\overline{l+1}$ does not appear in $\tilde{x}$ and $\overline{l}$ appears in $\tilde{x}$, which means that the properties in (d) hold. Hence $\rho(\tilde{x})=x$. The arguments for the reverse direction are the same .
\endproof
\end{lem}
Let \begin{equation}\label{1}f_l(\mathbf x):=\max\{1\leq p\leq k|\sigma^p_l(\mathbf x)-\tau^p_l(\mathbf x)=\min\{\sigma^j_l(\mathbf x)-\tau^j_l(\mathbf x)|1\leq j\leq k\}\}\end{equation}
\begin{equation}\label{2}e_l(\mathbf x):=\min\{1\leq p\leq k|\sigma^p_l(\mathbf x)-\tau^p_l(\mathbf x)=\min\{\sigma^j_l(\mathbf x)-\tau^j_l(\mathbf x)|1\leq j\leq k\}\}.\end{equation}
Now we are able to define the Kashiwara operators,
\begin{equation}\label{kashif}\tilde{f}_l\mathbf x=\begin{cases} 0,& \text{if $\theta_l(x_{q_{f_l(\mathbf x)}})=0$}\\
(\cdots,\theta_l(x_{q_{f_l(\mathbf x)}}),\cdots),& \text{else}\end{cases}\end{equation}
\begin{equation}\label{kashie}\tilde{e}_l\mathbf x=\begin{cases}0,& \text{ if $\rho_l(x_{q_{e_l(\mathbf x)}})=0$}\\
(\cdots,\rho_l(x_{q_{e_l(\mathbf x)}}),\cdots).& \text{else}\end{cases}\end{equation}
Let us consider an example:
\begin{ex}\mbox{}
\begin{enumerate}
\item Let $\Lg=A_4$ and $\mathbf x=(x_1,x_2,x_3,0,0,\cdots)$ with $$x_1=(2,2)(1,3)\in\mathcal{R}^{2}_2,\ x_2=(2,4)\in\mathcal{R}^{4}_1,\ x_3=(3,3)(2,4)\in\mathcal{R}^{3}_2.$$
For $l=2$ we get $\sigma^1_2(\mathbf x)=\sigma^2_2(\mathbf x)=\sigma^3_2(\mathbf x)=\tau^1_2(\mathbf x)=0,\tau^2_2(\mathbf x)=1,\tau^3_2(\mathbf x)=2$ and hence $$\tilde{f}_2\mathbf x=0.$$
\item Let $\Lg=C_3$ and $\mathbf x=(x_1,x_2,0,0,\cdots)$ with $$x_1=(1,\overline{2})\in\mathcal{R}^{2}_1,\ x_2=(3,3)(1,\overline{2})\in\mathcal{R}^{3}_2,$$ then $\sigma^2_3(\mathbf x)=0,\tau^2_3(\mathbf x)=1$ and $$\tilde{e}_3\mathbf x=(x_1,x_1,0,0\cdots).$$
\end{enumerate}
\end{ex}
It remains to define the maps $\varphi_l$ and $\epsilon_l$. These maps are given by the next formula:
\begin{equation}\label{epsilon}\epsilon_l(\mathbf x)=\tau^2_l((\emptyset_1,x_{q_1},0,0\cdots))-\min\{\sigma^j_l(\mathbf x)-\tau^j_l(\mathbf x)|1\leq j\leq k\}\end{equation}
\begin{equation}\label{phi}\varphi_l(\mathbf x)=\sigma^{k+1}_l(\mathbf x)-\tau^k_l(\mathbf x)-\min\{\sigma^j_l(\mathbf x)-\tau^j_l(\mathbf x)|1\leq j\leq k\}.\end{equation}
If we collect all the maps stated in (\ref{weight}), (\ref{kashif}), (\ref{kashie}), (\ref{epsilon}) and (\ref{phi}) we can formulate the next proposition.
\begin{prop}\label{crystalprop}
The set $\mathcal{R}^{\infty}$ becomes a semiregular crystal.
\proof
In Lemma~\ref{semiregularity} we proof the semiregularity of $\mathcal{R}^{\infty}$, which ensures that (4) and (5) from Definition~\ref{abcrystal} hold. So, to verify the proposition, it is sufficient to prove (1),(2),(3) and (6), where (2) and (3) are easily checked with the help of Remark~\ref{remark1}. Let us start by proving (1); so let $\mathbf x\in\mathcal{R}^{\infty}$ be arbitrary with finitely many non-zero components, say $x_{q_1}\in\mathcal{R}_{s_1}^{j_1},\cdots,x_{q_k}\in \mathcal{R}_{s_k}^{j_k}$. Then we order these components in a way such that the first components are contained in $\bigcup_s\bigcup_{m<l}\mathcal{R}_s^{m}$ followed by components in $\bigcup_s\bigcup_{m>l}\mathcal{R}_s^{m}$ and the last ones are contained in $\bigcup_s\mathcal{R}_s^{l}$. So we can write the set of non-zero components of $\mathbf x$ as a disjoint union of three subsets $A_{<l}\cup A_{>l}\cup A_{=l}$. Further let $\mathbf x_1$ be the element in $\mathcal{R}^{\infty}$ obtained from $\mathbf x$ by replacing all components not belonging to $A_{<l}$ by $0$. And $\mathbf x_2$ and $\mathbf x_3$ respectively are similar defined using $A_{>l}$ and $A_{=l}$, respectively. Then we get 
\begin{eqnarray*}&\varphi_l(\mathbf x)-\epsilon_l(\mathbf x)=\sigma^{k+1}_l(\mathbf x)-\tau^k_l(\mathbf x)-\tau^2_l((\emptyset_1,x_{q_1},0,0\cdots))&\\&=\sigma^{|A_{<l}|+1}_l(\mathbf x_1)-\tau^{|A_{<l}|}_l(\mathbf x_1)+\sigma^{|A_{>l}|+1}_l(\mathbf x_2)-\tau^{|A_{>l}|}_l(\mathbf x_2)&\\&+\sigma^{|A_{=l}|+1}_l(\mathbf x_3)-\tau^{|A_{=l}|}_l(\mathbf x_3)-\sum^3_{i=1}\tau^2_l((\emptyset_1,y_i,0,0\cdots)),\end{eqnarray*}
where the sequence $y_i$ is the first non-zero element in $\mathbf x_i$.
A short calculation shows 
\begin{align*}&\sigma^{|A_{<l}|+1}_l(\mathbf x_1)-\tau^{|A_{<l}|}_l(\mathbf x_1)-\tau^2_l((\emptyset_1,y_1,0,0\cdots))=-\sum_{1\leq i\leq k : x_{q_i}\in A_{<l}} \langle\alpha_l^\vee,\wt(x_{q_i})\rangle&\\&
\sigma^{|A_{>l}|+1}_l(\mathbf x_2)-\tau^{|A_{>l}|}_l(\mathbf x_2)-\tau^2_l((\emptyset_1,y_2,0,0\cdots))=-\sum_{1\leq i\leq k:x_{q_i}\in A_{>l}} \langle\alpha_l^\vee,\wt(x_{q_i})\rangle
&\\&\sigma^{|A_{=l}|+1}_l(\mathbf x_3)-\tau^{|A_{=l}|}_l(\mathbf x_3)-\tau^2_l((\emptyset_1,y_3,0,0\cdots))=c_l-\sum_{1\leq i\leq k:x_{q_i}\in A_{=l}} \langle\alpha_l^\vee,\wt(x_{q_i})\rangle. \end{align*}
Now we proceed to prove (6). By the definition of the Kashiwara operators and Lemma~\ref{vor1} it is enough to show $f_l(\tilde{e}_l\mathbf x)=e_l(\mathbf x)$ and $e_l(\tilde{f}_l\mathbf x)=f_l(\mathbf x)$. Since the proofs are similar we prove only the latter equation. Assume that $\tilde{f}_l\mathbf x\neq 0$ and $f_l(\mathbf x)\neq 1$, then $$\sigma^j_l(\tilde{f}_l\mathbf x)=\begin{cases}\sigma^j_l(\mathbf x),& \text{ if $j\leq f_l(\mathbf x)$}\\
\sigma^j_l(\mathbf x)-1,& \text{ if $j>f_l(\mathbf x)$}\end{cases},\quad \tau^j_l(\tilde{f}_l\mathbf x)=\begin{cases}\tau^j_l(\mathbf x),& \text{ if $j\leq f_l(\mathbf x)-1$}\\
\tau^j_l(\mathbf x)+1,& \text{ if $j\geq f_l(\mathbf x)$}.\end{cases}$$
Subsequently we have 
\begin{flalign*}&e_l(\tilde{f}_l\mathbf x)=\min\{1\leq p\leq k|\sigma^p_l(\tilde{f}_l\mathbf x)-\tau^p_l(\tilde{f}_l\mathbf x)=\min\{\sigma^j_l(\tilde{f}_l\mathbf x)-\tau^j_l(\tilde{f}_l\mathbf x)|1\leq j\leq k\}\}
&\\&=\min\{1\leq p\leq k|\sigma^p_l(\tilde{f}_l\mathbf x)-\tau^p_l(\tilde{f}_l\mathbf x)=\sigma^{f_l(\mathbf x)}_l(\mathbf x)-\tau^{f_l(\mathbf x)}_l(\mathbf x)-1\}=f_l(\mathbf x).\end{flalign*}
In the case where $f_l(\mathbf x)=1$, we obtain 
$$e_l(\tilde{f}_l\mathbf x)=\min\{1\leq p\leq k|\sigma^p_l(\tilde{f}_l\mathbf x)-\tau^p_l(\tilde{f}_l\mathbf x)=0\}=1.$$
\endproof
\end{prop}
Hence we have shown that the set $\mathcal{R}^{\infty}$ is an abstract crystal provided the semiregularity is shown. Thus our aim now is to verify that the maps $\varphi_l$ and $\epsilon_l$ respectively determine how often one can act with $\tilde{f}_l$ and $\tilde{e}_l$, respectively. The semiregularity is a necessary condition of a crystal $B$, if one wants to identify it with the crystal graph $B(\lambda)$.
\begin{lem}\label{semiregularity}
Let $\epsilon_l$ and $\varphi_l$ as in (\ref{epsilon}) and (\ref{phi}). For a given element $\mathbf x\in \mathcal{R}^{\infty}$ we obtain $$\epsilon_l(\mathbf x)=\max\{k\geq0|\tilde{e}_{l}^{k}\mathbf x\neq 0\},\quad \varphi_l(\mathbf x)=\max\{k\geq0|\tilde{f}_{l}^{k}\mathbf x\neq 0\}.$$
\proof
We proof the statement by induction on $z:=\max\{k\geq0|\tilde{e}_{l}^{k}\mathbf x\neq 0\}$; so let $z=0$ and suppose first that $e_l(\mathbf x)\neq 1$. By the definition of $e_l(\mathbf x)$ and $\rho_l$ we have $$\sigma_l^{e_l(\mathbf x)}(\mathbf x)-\tau^{e_l(\mathbf x)}_l(\mathbf x)<\sigma_l^{e_l(\mathbf x)-1}(\mathbf x)-\tau^{e_l(\mathbf x)-1}_l(\mathbf x)\Longleftrightarrow \sigma_l^{e_l(\mathbf x)}(\mathbf x)<\sigma_l^{e_l(\mathbf x)-1}(\mathbf x),$$ which is a contradiction to $e_l(\mathbf x)\neq 1$. Hence we have $e_l(\mathbf x)=1$ and as a consequence we obtain $\epsilon_l(\mathbf x)=\tau^2_l((\emptyset_1,x_{q_1},0,0\cdots))=0$, which proves the initial step. Now assume that $z>0$ and consider the element $\tilde{e}_l\mathbf x$, where we again presume initially $e_l(\mathbf x)\neq 1$. By applying the induction hypothesis and using $$\sigma^j_l(\tilde{e}_l\mathbf x)=\begin{cases}\sigma^j_l(\mathbf x),& \text{ if $j\leq e_l(\mathbf x)$}\\
\sigma^j_l(\mathbf x)+1,& \text{ if $j>e_l(\mathbf x)$}\end{cases},\quad \tau^j_l(\tilde{e}_l\mathbf x)=\begin{cases}\tau^j_l(\mathbf x),& \text{ if $j\leq e_l(\mathbf x)-1$}\\
\tau^j_l(\mathbf x)-1,& \text{ if $j\geq e_l(\mathbf x)$}\end{cases}$$ we arrive at \begin{flalign*}\max\{k\geq0|\tilde{e}_{l}^{k+1}\mathbf x\neq 0\}=\epsilon_l(\tilde{e}_l\mathbf x)=\tau^2_l((\emptyset_1,x_{q_1},0,0\cdots))-(\sigma^{e(\mathbf x)}_l(\mathbf x)-\tau^{e(\mathbf x)}_l(\mathbf x)+1)=\epsilon_l(\mathbf x)-1.\end{flalign*}
If $e_l(\mathbf x)= 1$, $$\epsilon_l(\tilde{e}_l(\mathbf x))=\tau^2_l((\emptyset_1,\rho_l(x_{q_1}),0,0\cdots))+\sigma^{e(\mathbf x)}_l(\mathbf x)-\tau^{e(\mathbf x)}_l(\mathbf x)=\epsilon_l(\mathbf x)-1.$$
The proof of the remaining equality $\varphi_l(\mathbf x)=\max\{k\geq0|\tilde{f}_{l}^{k}\mathbf x\neq 0\}$ is quite parallel.
\endproof
\end{lem}
Before we introduce the subcrystals $\mathcal{R(\lambda)}$ we need some facts about the theory of tensor products of crystals. The tensor product rule is a very nice combinatorial feature and important to realize the crystal bases of a tensor product of two $\U_q(\Lg)$-modules. 
%
\section{Tensor products and Nakajima monomials}\label{section4}
In this section, we want to recall tensor products of crystals and investigate the action of Kashiwara operators on tensor products. With the aim to have a different realization of $B(\lambda)$ from our approach we want to introduce the set of all Nakajima monomials, such that we can think of $B(\lambda)$ in terms of certain monomials. This theory is discovered by Nakajima \cite{Nak03}, and generalized by Kashiwara \cite{Kashi03}.
%
\subsection{Tensor product of crystals}
Suppose that we have two abstract crystals $B_1$, $B_2$ in the sense of Definition~\ref{abcrystal}, then we can construct a new crystal which is as a set $B_1\times B_2$. This crystal is denoted by $B_1\otimes B_2$ and the Kashiwara operators are given as follows:
$$\tilde{f}_{l}(b_1\otimes b_2)=\begin{cases} (\tilde{f}_{l}b_1)\otimes b_2, \text{ if $\varphi_l(b_1)>\epsilon_l(b_2)$}\\
b_1\otimes (\tilde{f}_{l}b_2), \text{ if $\varphi_l(b_1)\leq \epsilon_l(b_2)$}\end{cases}$$
$$\tilde{e}_{l}(b_1\otimes b_2)=\begin{cases} (\tilde{e}_{l}b_1)\otimes b_2, \text{ if $\varphi_l(b_1)\geq \epsilon_l(b_2)$}\\
b_1\otimes (\tilde{e}_{l}b_2), \text{ if $\varphi_l(b_1)<\epsilon_l(b_2)$.}\end{cases}$$
Furthermore, one can describe explicitly the maps $\wt$, $\varphi_l$ and $\epsilon_l$ on $B_1\otimes B_2$, namely:
$$\wt(b_1\otimes b_2)=\wt(b_1)+\wt(b_2)$$
$$\varphi_l(b_1\otimes b_2)=\max\{\varphi_l(b_2),\varphi_l(b_1)+\varphi_l(b_2)-\epsilon_l(b_2)\}$$
$$\epsilon_l(b_1\otimes b_2)=\max\{\epsilon_l(b_1),\epsilon_l(b_1)+\epsilon_l(b_2)-\varphi_l(b_1)\}.$$
One of the most important interpretation of the tensor product rule is the following theorem (for more details see \cite{HK02}).
\begin{thm}\label{}
Let $M_j$ be a integrable module in the category $\mathcal{O}^q$ and let $(L_j,B_j)$ be a crystal bases of $M_j, (j=1,2)$. Set $L=L_1\otimes L_2$ and $B=B_1\otimes B_2$. Then $(L,B)$ is a crystal bases of $M_1\otimes M_2$.
\end{thm}
%
\subsection{Nakajima monomials}
For $i\in I$ and $n\in\Z$ we consider monomials in the variables $Y_i(n)$, i.e. we obtain the set of Nakajima monomials $\mathcal{M}$ as follows:
$$\mathcal{M}:=\{\prod_{i\in I, n\in\Z} Y_i(n)^{y_i(n)}| y_i(n)\in\Z \mbox{ vanish except for finitely many $(i,n)$}\}$$
With the goal to define a crystal structure on $\mathcal{M}$, we take some integers $c=(c_{i,j})_{i\neq j}$ such that $c_{i,j}+c_{j,i}=1$. Let now $M=\prod_{i\in I, n\in\Z} Y_i(n)^{y_i(n)}$ be an arbitrary monomial in $\mathcal{M}$ and $l\in I$, then we set:
$$\wt(M)=\sum_i(\sum_n y_i(n))\omega_i$$
$$\varphi_l(M)=\max\{\sum_{k\leq n}y_l(k)|n\in\Z\},\quad \epsilon_l(M)=\max\{-\sum_{k>n}y_l(k)|n\in\Z\}$$
and 
$$n^l_f=\min\{n|\varphi_l(M)=\sum_{k\leq n}y_l(k)\},\quad n^l_e=\max\{n|\epsilon_l(M)=-\sum_{k>n}y_l(k)\}.$$
The Kashiwara operators are defined as follows:
$$\tilde{f}_{l}M=\begin{cases}
A_l(n^l_f)^{-1}M, \text{ if $\varphi_l(M)>0$}\\
0, \text{ if $\varphi_l(M)=0$}\end{cases},\quad \tilde{e}_{l}M=\begin{cases}
A_l(n^l_e)M, \text{ if $\epsilon_l(M)>0$}\\
0, \text{ if $\epsilon_l(M)=0$,}\end{cases}$$
whereby $$A_l(n):=Y_l(n)Y_l(n+1)\prod_{i\neq l}Y_i(n+c_{i,l})^{\langle \alpha_i^{\vee},\alpha_l\rangle}.$$
The following two results are shown by Kashiwara \cite{Kashi03}:
\begin{prop}
With the maps $\wt$, $\varphi_l$, $\epsilon_l$, $\tilde{f}_{l}$, $\tilde{e}_{l}$, $l\in I$, the set $\mathcal{M}$ becomes a semiregular crystal.
\end{prop}

\begin{rem}
\textit{A priori} the crystal structure depends on $c$, hence we will denote this crystal by $\mathcal{M}_c$. But it is easy to see that the isomorphism class of $\mathcal{M}_c$ does not depend on this choice. In the literature $c$ is often chosen as
$$c_{i,j}=\begin{cases} 0, \text{ if $i>j$}\\
 1, \text{ else} \end{cases} \quad or \quad c_{i,j}=\begin{cases} 0, \text{ if $i<j$}\\
 1, \text{ else.} \end{cases}$$
\end{rem}
\begin{prop}\label{iii}
Let $M$ be a monomial in $\mathcal{M}$, such that $\tilde{e}_{l}M=0$ for all $l\in I$. Then the connected component of $\mathcal{M}$ containing $M$ is isomorphic to $B(\wt(M))$.
\end{prop}
According to the latter proposition, it is of great interest to describe these connected components explicitly. This is worked out for special highest weight monomials for all classical Lie algebras in \cite{KKS2004},\cite{KKS04} and for the affine Lie algebra $A^{(1)}_n$ in \cite{K05}. We recall the results here only for type $C_n$ stated originally in \cite{KKS2004}.
\begin{prop}\label{iiii}
Let $\lambda=\sum^n_{i=1}m_i\omega_i$ be a dominant integral weight and consider the highest weight monomial $M=Y_1(1)^{m_1}\cdots Y_n(1)^{m_n}$. Then the connected component of $\mathcal{M}$ containing $M$ is characterized as the set of monomials of the form $$X_{t_{1,1}}(1)\cdots X_{t_{1,\alpha_1}}(1)\cdots X_{t_{n,1}}(n)\cdots X_{t_{n,\alpha_n}}(n),\quad t_{r,s}\in\{1<\cdots<n<\overline{n}<\cdots<\overline{1}\},$$ satisfying
\begin{enumerate}
\item $\alpha_j=m_{j}+\cdots+m_n$ for all $j=1,\ldots,n,$
\item $t_{j,1}\geq\ldots\geq t_{j,\alpha_j}$ for all $j=1,\ldots,n,$
\item $t_{j-1,k}>t_{j,k}$ for all $j=2,\ldots,n$ and $k=1,\ldots,\alpha_j,$ where 
$$X_i(m)=Y_{i-1}(m+1)^{-1}Y_i(m),\quad X_{\overline{i}}(m)=Y_{i-1}(m+(n-i+1))Y_i(m+(n-i+1))^{-1}.$$
\end{enumerate}
\end{prop}
Summarized, we have a semiregular crystal $\mathcal{M}$ and for each dominant integral weight $\lambda$ certain connected subcrystals contained in $\mathcal{M}$. These are isomorphic to $B(\lambda)$ and an explicit description of these components is worked out for the classical simple Lie algebras and $A^{(1)}_n$.
In the remaining parts of this paper we prove a similar result as Proposition~\ref{iii} and Proposition~\ref{iiii}, whereby our ``big" semiregular crystal is $\mathcal{R}^{\infty}$.
%
\section{Explicit description of the connected components}\label{section5}
In this section we define for the dominant integral weights $\lambda=\sum^n_{i=1}m_i\omega_i$ certain connected subcrystals $\mathcal{R}(\lambda)\subseteq \mathcal{R}^{\infty}$. Furthermore, we provide an explicit description of these crystals in Theorem~\ref{char} and Theorem~\ref{char2} respectively, i.e. we give a set of conditions describing $\mathcal{R}(\lambda)$.

\begin{defn}\label{rlambda}
For a dominant integral weight $\lambda=\sum^n_{i=1}m_i\omega_i$ let $\mathcal{R}(\lambda)$ be the connected component of $\mathcal{R}^{\infty}$ containing $$r_{\lambda}=(\underbrace{\emptyset_1,\cdots,\emptyset_1}_{m_1},\cdots,\underbrace{\emptyset_n,\cdots,\emptyset_n}_{m_n},0,0\cdots).$$
\end{defn}
Note that the weight of $r_\lambda$ is precisely $\lambda$. Furthermore, by definition, $\mathcal{R}(\lambda)$ is connected and for $\lambda=\omega_i$ we can immediately provide a description of $\mathcal{R}(\omega_i)$. To be more accurate we prove as a first step the following proposition:
\begin{prop}\label{prop1}
Let $\lambda=\omega_i$, then we obtain $$\mathcal{R}(\lambda)=(\bigcup_s R^i_s,0,0,\cdots)=\bigcup_s R^i_s.$$
\proof
Since $(\bigcup_s R^i_s,0,0,\cdots)$ is stable under the Kashiwara operators $\tilde{f}_l$ and $\tilde{e}_l$ (Remark~\ref{remark1}) it is enough to prove that $\emptyset_i$ is the unique highest weight element in $\bigcup_s R^i_s$, i.e. $$\tilde{e}_l(x,0,0,\cdots)=0 \ \forall l\in I \Rightarrow x=(\emptyset_i,0,0\cdots).$$ 
Assume $x=(i_1,i_1')\cdots (i_s,i'_s)$ and $\tilde{e}_l(x,0,0,\cdots)=0 \ \forall l\in I$. If $i_1< i$ we have $\tilde{e}_{i_1}(x,0,0,\cdots)\neq 0$ and if $i < i'_1\leq n$ we have $\tilde{e}_{i'_1}(x,0,0,\cdots)\neq 0$. So let $\overline{n-1}\leq i'_1$ (this case cannot appear if $\Lg=A_n$), then since $i'_1\in\{i'_1,\cdots,i'_s\}$, $\overline{\overline{i'_1}+1}\notin\{i'_1,\cdots,i'_s\}$ and $i'_1\leq\overline{i_1}=\overline{i}$ we obtain that $\overline{i'_1}\in\{i,\cdots,n-1\}$. Thus $x$ satisfies (d) and hence $\tilde{e}_{\overline{i'_1}}(x,0,0,\cdots)\neq 0$. According to these calculations, the pair $(i,i)$ must appear in $x$. This implies $\tilde{e}_{i}(x,0,0,\cdots)\neq 0$.
\endproof
\end{prop}
For general $\lambda$'s we can describe the connected components and refer to the following two subsections.
%
\subsection{Explicit description of \texorpdfstring{$\mathcal{R}(\lambda)$}{R} in type \texorpdfstring{$A_n$}{A}}
In this subsection we give an explicit characterization of $\mathcal{R}(\lambda)$ if $\Lg$ is the special linear Lie algebra, i.e. we give conditions whether a sequence $\mathbf x$ is contained in $\mathcal{R}(\lambda)$ or not (Theorem~\ref{char}). Initially we note that $\mathcal{R}(\lambda)$ lives in $\mathcal{R}^{k=\sum_{i}m_i}$ and the first $k$ components are non-zero. For simplicity we set $e_l(\mathbf x)=e(\mathbf x)$, $f_l(\mathbf x)=f(\mathbf x)$, $\theta_l=\theta$ and $\rho_l=\rho$.

\begin{thm}\label{char}
The crystal $\mathcal{R}(\lambda)$ consists of all sequences $$\mathbf x=(x_1,\cdots,x_k,0,0,\cdots),$$ such that:
\begin{enumerate}
\item for all pairs $(x_q\neq \emptyset_i,x_{q+1}\neq \emptyset_j)$, say $x_q=(i_1,i'_1)\cdots(i_s,i'_s)\in \mathcal{R}^i_s$ and $x_{q+1}=(j_1,j'_1)\cdots(j_t,j'_t)\in\mathcal{R}^j_t$ we have
\begin{enumerate}[i)]
\item $\sharp\{i_k|1\leq i_k\leq j_{t-p}\}\geq p+1$, for all $0\leq p\leq u:=\max \{0\leq r\leq t-1|j_{t-r}\leq i\}.$
\item $j_{t-p}\leq i'_{j_{t-p}-p+s-i}$, for all $p\in\{u+1\leq r\leq t-1|j_{t-r}\leq r+i\}$
\item $j'_p\leq i'_{j-i+s-t+p}$, for all $1\leq p \leq i-j+t$,
\end{enumerate}
\item there is no pair $(x_q,x_{q+1})$ of the form $(\emptyset_i,x_{q+1}\neq \emptyset_j)$ with $i\geq j_t$.
\end{enumerate}
\proof
First we note that the element $r_{\lambda}$ is contained in $\mathcal{R}(\lambda)$ and is a highest weight element. Furthermore we claim that $r_{\lambda}$ is the unique highest weight element. So suppose that we have another element $\mathbf x=(x_1,\cdots,x_k,0,0,\cdots)$ satisfying $\tilde{e}_l\mathbf x=0$ for all $l\in I$. Let $z$ be the lowest integer which appears in one of the sequences $x_1,\cdots,x_k$ and let $p$ be the minimal integer such that $z$ appears in $x_p$, say $x_p=(j_1,j'_1)\cdots(j_t=z,j'_t)\in \mathcal{R}^j_t$. In the case where such a $z$ does not exist we have $\mathbf x=r_{\lambda}$. We remark that \begin{equation}\label{zw}\tilde{e}_l(x_p,0,0\cdots)=0\mbox{ for all $l\geq z$}\end{equation} would imply the claim, whereby the reason is the following:\par
assume that (\ref{zw}) holds and let $r=\max\{1\leq r \leq t|j_{r-1}\neq j_{r}+1\}$. If $j_1=j$ and $r=1$ we set $i=j'_1+1$ and $i=j_r+1$ else. In either case we obtain $\epsilon_{i-1}(x_p,0,0,\cdots)\neq 0$, which is a contradicition to (\ref{zw}). So if suffices to show (\ref{zw}).\par
Primarily we claim that $x_1,\cdots,x_{p-1}\in\{\emptyset_1,\cdots,\emptyset_n\}$. Assume that the element $x_{p-1}=(i_1,i'_1)\cdots(i_s,i'_s)\in\mathcal{R}^i_s$ is not in the aforementioned set; then the property $(1)i)$ is violated by the choice of $z$ and $p$, because $$0=\sharp\{i_k|1\leq i_k\leq z\}\geq 1.$$ Hence, there exists a pair $(\emptyset_j,x_p)$ which forces $j< z$ and thus, again by the choice of $z$, the remaining sequences are again contained in the set $\{\emptyset_1,\cdots,\emptyset_n\}.$
As a consequence, a short calculation by using Lemma~\ref{semiregularity} shows the uniqueness of the highest weight element, namely $$\epsilon_l(\mathbf x)=0 \ \forall l\in I \Rightarrow \epsilon_l((x_p,\cdots,x_k,0,0\cdots))=0 \ \forall l\geq z \Rightarrow \epsilon_l((x_p,0,0\cdots))=0\mbox{ for all $l\geq z$}.$$ 
In order to obtain a connected crystal it remains to show that $\tilde{e}_l\mathbf x, \tilde{f}_l\mathbf x\in\mathcal{R}(\lambda)\cup \{0\}.$ Assume that $\tilde{f}_l\mathbf x\neq 0$, say $$\tilde{f}_l\mathbf x=(\cdots,x_{f(\mathbf x)-1},\theta(x_{f(\mathbf x)}),x_{f(\mathbf x)+1},\cdots),$$ where we set for simplicity $f(\mathbf x)=q$. Our goal here is to show that the properties (1) and (2) hold for $\tilde{f}_l\mathbf x$, where we start by proving (2).\par
It is easy to see that (2) can only be violated if one of the following two cases occur
\begin{align*}
&\bullet\ (x_{q-1},x_{q})=(\emptyset_i,\emptyset_i) \mbox{ and } l=i
\notag&\\& \bullet\ x_q=(i_1,i'_1)\cdots(i_s,i'_s)\in \mathcal{R}^i_s, x_{q-1}=\emptyset_{i_s-1} \mbox{ and } l=i_s-1<i.
\end{align*}
In either case we obtain $\sigma_l^{q-1}(\mathbf x)-\tau^{q-1}_l(\mathbf x)<\sigma^q_l(\mathbf x)-\tau^q_l(\mathbf x)$, which is a contradiction to the choice of $q$. The proof of the fact that (1) holds will proceed in several cases. In the remaining parts we denote the entries of $\theta(x_q)$ by $\theta(i_k),\theta(i'_k)$, i.e. 
$$\theta(x_q)=(\theta(i_1),\theta(i'_1))\cdots(\theta(i_{s+\delta_{l,i}}),\theta(i'_{s+\delta_{l,i}})).$$\vskip 5pt
\textbf{Case 1.1:}\ (1) is violated for the pair $(\theta(x_q),x_{q+1})$:\\
Let $x_q=(i_1,i'_1)\cdots(i_s,i'_s)\in \mathcal{R}^i_s$ and $x_{q+1}=(j_1,j'_1)\cdots(j_t,j'_t)\in\mathcal{R}^j_t$.
We first consider the case $l<i$, which means that we replace $l+1$ by $l$. Since the entries $i'_p$ stay unchanged, only property (1)i) can be violated. However, property (1)i) is still fulfilled, because $$\sharp\{i_k|1\leq i_k\leq j_{t-p}\}\leq \sharp\{\theta(i_k)|1\leq \theta(i_k)\leq j_{t-p}\}, \mbox{ for all } 0\leq p\leq u.$$\par
If $l>i$ we replace $l-1$ by $l$ and hence (1) is obviously fulfilled. So suppose that $l=i$, which means that we add the entry $(l,l)$. The equality $$\sharp\{i_k|1\leq i_k\leq j_{t-p}\}=\sharp\{\theta(i_k)|1\leq \theta(i_k)\leq j_{t-p}\},\quad \forall 0\leq p\leq u-1$$ and inequality $$\sharp\{i_k|1\leq i_k\leq j_{t-u}\}\leq\sharp\{\theta(i_k)|1\leq \theta(i_k)\leq j_{t-u}\}$$ imply that property (1)i) is still fulfilled. Furthermore, since $\theta(i'_k)=i'_{k-1}$ for $k>1$, $\theta(i'_1)=i$ and $\theta(x_q)\in\mathcal{R}^i_{s+1}$, we obtain 
\begin{align*}&j_{t-p}\leq i'_{j_{t-p}-p+s-i}=\theta(i'_{j_{t-p}-p+s+1-i}),\mbox{ for all } p\in\{u+1\leq r\leq t-1|j_{t-r}\leq r+i\}&\\&
j'_p\leq i'_{j-i+s-t+p}=\theta(i'_{j-i+s+1-t+p}), \mbox{ for all }  1\leq p \leq i-j+t,\end{align*}which particularly means that Case 1.1 can never appear.\vskip3pt
\textbf{Case 1.2:}\ (1)i) is violated for the pair $(x_{q-1},\theta(x_q))$:\\ 
For simplicity we set $x_{q-1}=(i_1,i'_1)\cdots(i_s,i'_s)\in \mathcal{R}^i_s$ and $x_{q}=(j_1,j'_1)\cdots(j_t,j'_t)\in\mathcal{R}^j_t$. In that case there exists at least one $0\leq p\leq \theta(u)=\max \{r|\theta(j_{t-r})\leq i\}$, such that \begin{equation}\label{vio}\sharp\{i_k|1\leq i_k\leq \theta(j_{t-p})\}<p+1.\end{equation} Necessarily we must have $l\leq i$ and in the case where $l<j$ (\ref{vio}) implies \begin{eqnarray*}&\sharp\{i_k|1\leq i_k\leq l+1)\}>\sharp\{i_k|1\leq i_k\leq l)\}\notag&\\&\ \sharp\{i_k|1\leq i_k\leq l-1)\}=\sharp\{i_k|1\leq i_k\leq l)\}.\end{eqnarray*} As a consequence we get that $l+1$ appears in $x_{q-1}$ while $l$ does not appear and thus $\sigma_l^{q-1}(\mathbf x)-\tau^{q-1}_l(\mathbf x)<\sigma^q_l(\mathbf x)-\tau^q_l(\mathbf x)$, which is a contradiction to the choice of $q$.\par
Eventually if $l=j=i$ we obtain in a similar way \begin{equation}t\geq \sharp\{i_k|1\leq i_k\leq \theta(j_{1})=j\}=s\geq\sharp\{i_k|1\leq i_k\leq j_1)\}\geq t,\end{equation}which forces on the one hand $s=t$ and on the other hand that $j$ does not appear in $x_{q-1}$. To be more precise, we can conclude the latter statement with the help of (1)iii) and (1)i), namely $$j< j'_1\leq i'_{s-t+1}=i'_1$$ and
$$s\geq \sharp\{i_k|1\leq i_k< i)\}\geq \sharp\{i_k|1\leq i_k\leq j_1)\}\geq t=s.$$ Thus we get again $\sigma_l^{q-1}(\mathbf x)-\tau^{q-1}_l(\mathbf x)<\sigma^q_l(\mathbf x)-\tau^q_l(\mathbf x)$, which is once more a contradiction to the choice of $q$.\vskip3pt
\textbf{Case 1.3:}\ (1)ii) is violated for the pair $(x_{q-1},\theta(x_q))$:\\
Here we have $l\leq j$ and if $l=j$ we must have $j>i$, because otherwise (1)ii) wouldn't be violated. We first consider the case where $l=j$ and notice, that the only possible violation is given by the following inequality $$j>i'_{j-t+s-i}.$$ We can conclude that $j$ does not occur in $x_{q-1}$, because either $j=t+i$ and hence $$i'_s<j
$$ or $j+1\leq t+i$ and (1)iii) is applicable, which yields $j<j'_1\leq i'_{j-t+s-i+1}$. To obtain a contradiction we have to show that $j-1$ appears in $x_{q-1}$. If $j_1>i$, this follows by the subsequent calculation:
$$i'_{j-t+s-i}\leq j-1 \leq j_1+(j-1-j_1)\leq i'_{j_1-t+1+s-i}+(j-1-j_1)\leq i'_{j-t+s-i}.$$
If $j_1\leq i$ we obtain with property (1)i) that $s\geq\{i_k|1\leq i_k\leq j_1\}\geq t$. In particular we actually have $s=t$, because otherwise we would get $$j\leq j-t+s-1\leq i'_{j-t+s-i}.$$ Eventually we can conculde again that $j-1$ must appear in $x_{q-1}$ $$i'_{j-i}\leq j-1\leq i'_{j-i}.$$
Now we suppose $l<j$. Then an easy consideration shows that (1)ii) can only be violated if $l>i$ and thus we obtain similar as before that the only violation which can occur is the following $$l>i'_{l-r+s-i},$$ where we expect $j_{t-r}=l+1.$
We would like to show as before that $l$ does not appear in $x_{q-1}$ while $l-1$ appears. We either have $l=r+i$ and thus $l>i'_s$ or $l+1\leq r+i$. In the latter case we apply property (1)ii) and obtain $$l<i'_{l-r+s-i+1}.$$ In either case we notice that $l$ does not occur in $x_{q-1}$.\par
In order to prove the remaining part we consider the element $j_{t-(r-1)}$ which is considerable, since
$$l-r+s-1\leq i'_{l-r+s-i}\leq l-1\Longrightarrow 1\leq s\leq r.$$ 
In the case where $j_{t-(r-1)}$ is greater than $i$ we obtain
$$i'_{l-r+s-i}\leq l-1\leq j_{t-(r-1)}+(l-1- j_{t-(r-1)})\leq i'_{j_{t-(r-1)}-r+1+s-i}+(l-1- j_{t-(r-1)})\leq i'_{l-r+s-i}$$ and otherwise by using property (1)i) we get $r=s$. Thus $$l-1\geq i'_{l-i}\geq l-1.$$
\vskip3pt
\textbf{Case 1.4:}\ (1)iii) is violated for the pair $(x_{q-1},\theta(x_q))$:\\
We suppose that (1)iii) is violated, which forces $l\geq j$. In the case where $l>j$ we have $$\exists p : j'_p=l-1=i'_{j-i+s-t+p}.$$ It follows that $l-1$ appears and $l$ does not appear in $x_{q-1}$, because $l=i'_{j-i+s-t+p+1}$ would imply $$l=i'_{j-i+s-t+p+1}\geq j'_{p+1}>j'_{p}=l-1 \Longrightarrow j'_{p+1}=l.$$ Consequently, $\sigma_l^{q-1}(\mathbf x)-\tau^{q-1}_l(\mathbf x)<\sigma^q_l(\mathbf x)-\tau^q_l(\mathbf x)$ and we obtain as usual a contradiction to the choice of $q$.\par
If the remaining case $l=j$ occurs, then one of the inequalities $$\theta(j'_p)\leq i'_{j-i+s-t-1+p},\ 1\leq p \leq i-j+t+1,$$ must be violated. Clearly, the only possibility is that $j\leq i'_{j-i+s-t}$ does not hold. If $t=j-i$ we get $$j>i'_s$$ and else we can assume $j-i+1\leq t$ so that (1)iii) is applicable, which yields $$i\leq i'_{j-i+s-t}<j<j'_1\leq i'_{j-i+s-t+1}.$$ Therefore, $j$ does not appear in $x_{q-1}$. In what follows, we finish our proof by showing that $j-1$ appears in $x_{q-1}$. If $j_1>i$ we can apply (1)ii) and get 
$$i'_{j-i+s-t}\leq j_1+(j-1-j_1)\leq i'_{j_{1}-t+1+s-i}+(j-1-j_1)\leq i'_{j-i+s-t}.$$ If $j_1\leq i$ we can verify with (1)i) that $s\geq t$, but the assumption $s-1\geq t$ yield in a contradiction, namely $$i'_{j-i+s-t}\geq j+s-t-1\geq j.$$ So $s=t$ and we obtain the required equality $j-1\geq i'_{j-i}\geq j-1$.\par The proof of $\tilde{e}_l\mathbf x\in\mathcal{R}(\lambda)\cup \{0\}$ is similar, which completes the proof.
\endproof
\end{thm}
%
\subsection{Explicit description of \texorpdfstring{$\mathcal{R}(\lambda)$}{R} in type \texorpdfstring{$C_n$}{C}}
In this subsection we would like to give an explicit characterization of $\mathcal{R}(\lambda)$ if $\Lg$ is a symplectic Lie algebra. In order to state the main theorem we fix some notation.\par
For an arbitrary subset $A\subseteq\mathbf I$, $n\in\{0,1\}$ and $y\in\mathbf I$ we set $$\delta_{y,A}=\begin{cases}1,&\text{if $y\in A$}\\
0,&\text{if $y\notin A$}\end{cases},\quad nA=\begin{cases}A,&\text{if $n=1$}\\
\emptyset,&\text{if $n=0$.}\end{cases}$$
The analogue result to Theorem~\ref{char} for type $C_n$ is the following:
\begin{thm}\label{char2}
The crystal $\mathcal{R}(\lambda)$ consists of all sequences $$\mathbf x=(x_1,\cdots,x_k,0,0,\cdots),$$ such that
\begin{enumerate}
\item for all pairs $(x_q\neq \emptyset_i,x_{q+1}\neq \emptyset_j)$, say $x_q=(i_1,i'_1)\cdots(i_s,i'_s)\in \mathcal{R}^i_s$ and $x_{q+1}=(j_1,j'_1)\cdots(j_t,j'_t)\in\mathcal{R}^j_t$ we have
\begin{enumerate}[i)]
\item $\sharp\{i_k|1\leq i_k\leq j_{t-p}\}\geq p+1$, for all $0\leq p\leq u:=\max \{0\leq r\leq t-1|j_{t-r}\leq i\}.$
\item $j_{t-p}\leq i'_{j_{t-p}-p+s-i}$, for all $p\in\{u+1\leq r\leq t-1|j_{t-r}\leq r+i\}$
\item $j'_p\leq i'_{j-i+s-t+p}$, for all $1\leq p \leq i-j+t$,
\end{enumerate}
\item there is no pair $(x_q,x_{q+1})$ of the form $(\emptyset_i,x_{q+1}\neq \emptyset_j)$ with $i\geq j_t$,
\item there is no pair $(x_q\neq\emptyset_i,x_{q+1})$ with the following property:\vskip2pt
there exists $i'_p\geq i'_r\geq n$ with
\begin{enumerate}[(a)] 
\item $\overline{i'_r}-\delta_{\overline{i'_r},\{i+1,\cdots,n\}}\notin (1-\delta_{\overline{i'_r},\{i+1,\cdots,n\}})\{i_1,\ldots,i_s\}\ \cup \ \delta_{\overline{i'_r},\{i+1,\cdots,n\}}(\mathbf I- \{i'_1,\ldots,i'_s\})$
\item $\overline{i'_p}-\delta_{\overline{i'_p},\{j+1,\cdots,n\}}\notin (1-\delta_{\overline{i'_p},\{j+1,\cdots,n\}})\{j_1,\ldots,j_t\}\ \cup \ \delta_{\overline{i'_p},\{j+1,\cdots,n\}}(\mathbf I- \{j'_1,\ldots,j'_t\})$
\item$(p-r)+\sharp\{j_k|j_k<\overline{i'_p}\}-\sharp\{i_k|i_k<\overline{i'_r}\}+\sharp\{i'_k|i'_k<\overline{i'_r}\}-\sharp\{j'_k|j'_k<\overline{i'_p}\}\geq \max\{0,\overline{i'_r}-i\}-\max\{0,\overline{i'_p}-j\},$
\end{enumerate}
\item there is no pair $(x_q\neq\emptyset_i,x_{q+1}\neq\emptyset_j)$ with the following property:\vskip2pt
there exists $i'_p\geq j'_r\geq n$ with
\begin{enumerate}[(a)] 
\item $\overline{i'_p}-\delta_{\overline{i'_p},\{j+1,\cdots,n\}}\notin (1-\delta_{\overline{i'_p},\{j+1,\cdots,n\}})\{j_1,\ldots,j_t\}\ \cup \ \delta_{\overline{i'_p},\{j+1,\cdots,n\}}(\mathbf I- \{j'_1,\ldots,j'_t\})$
\item $\overline{j'_r}-\delta_{\overline{j'_r},\{j+1,\cdots,n\}}\notin (1-\delta_{\overline{j'_r},\{j+1,\cdots,n\}})\{j_1,\ldots,j_t\}\ \cup\ \delta_{\overline{j'_r},\{j+1,\cdots,n\}}(\mathbf I- \{j'_1,\ldots,j'_t\})$
\item $(i-j)+(t-s)+(p-r)+\sharp\{j'_k|\overline{i'_p}\leq j'_k<\overline{j'_r}\}-\sharp\{j_k|\overline{i'_p}\leq j_k<\overline{j'_r}\}\geq \max\{0,\overline{j'_r}-j\}-\max\{0,\overline{i'_p}-j\}.$
\end{enumerate}
\end{enumerate}
\proof
First of all we note as in the $A_n$ case, that the element $r_{\lambda}$ is contained in $\mathcal{R}(\lambda)$ and is a highest weight element. In order to prove that $r_{\lambda}$ is the unique highest weight element in $\mathcal{R}(\lambda)$ we assume $\mathbf x=(x_1,\cdots,x_k,0,0,\cdots)$ to be another one. Let $z$ be the lowest integer which appears in one of the sequences $x_1,\cdots, x_k$ and let $p$ be the minimal integer such that $z$ appears in $x_p$, say $x_p=(i_1,i'_1)\cdots(i_s=z,i'_s)\in \mathcal{R}^i_s$. 
We can prove similar to Theorem~\ref{char} that the elements $x_1,\ldots,x_p$ are contained in $\{\emptyset_1,\cdots,\emptyset_n\}$. Accordingly we get once more with Lemma~\ref{semiregularity} (similar to the $A_n$ case)
$$\tilde{e}_l\mathbf x=0 \ \forall l\in I \Rightarrow \tilde{e}_l(\mathbf (x_p,\cdots,x_k,0,0\cdots))=0 \ \forall l\geq z \Rightarrow \tilde{e}_l(\mathbf (x_p,0,0\cdots))=0\mbox{ for all $l\geq z$}.$$
Our aim is again to prove the impossibility of $\tilde{e}_l(\mathbf (x_p,0,0\cdots))=0\mbox{ for all $l\geq z$}$. Let $r=\max\{1\leq r \leq s|i_{r-1}\neq i_{r}+1\}$. If $i_1=i$, $r=1$ and $i'_1\notin\{\overline{n-1},\ldots,\overline{1}\}$ we set $j=i'_1+1$ and if $i_1=i$,$r=1$ is not satisfied we set $j=i_r+1$ and obtain similar to Theorem~\ref{char} that $\epsilon_{j-1}(x_p,0,0,\cdots)\neq 0$. Thus, the only remaining case which can appear is when $i_1=i$, $r=1$ and $i'_1\in\{\overline{n-1},\ldots,\overline{1}\}$. In this particular case we set $j=\overline{i'_1}$ and claim $\epsilon_{j}(x_p,0,0,\cdots)\neq 0$. The latter claim is true, because on the one hand we have 
$$ i'_1\in\{i'_1,\cdots,i'_s\}, \overline{\overline{i'_1}+1}\notin\{i'_1,\cdots,i'_s\}$$ 
and on the other hand we got $\overline{j}=i'_1\leq \overline{i_1}=\overline{i}\Longrightarrow j\geq i$, which verfies the properties listed in (d). To be more precise, if $j>i$ we have $j-1\notin\{i'_1,\cdots,i'_s\}$ and if $j=i$ we have $j=i=i_1$.\par
In order to finish the theorem it remains to show that $\mathcal{R}(\lambda)$ is stable under the Kashiwara operators, i.e. $\tilde{e}_l\mathbf x, \tilde{f}_l\mathbf x\in\mathcal{R}(\lambda)\cup \{0\}.$ Assume that $\tilde{f}_l\mathbf x\neq 0$, say $$\tilde{f}_l\mathbf x=(\cdots,x_{q-1},\theta(x_{q}),x_{q+1},\cdots).$$ Our aim here is to prove that the properties (1)-(4) hold for $\tilde{f}_l\mathbf x$, whereby the verification of the first and second property proceeds almost similar to Theorem~\ref{char}. Nevertheless we will demonstrate some parts of it in Case 1. In the remaining parts of our proof we set $x_q=(i_1,i'_1)\cdots(i_s,i'_s)\in \mathcal{R}^i_s$ and $x_{q+1}=(j_1,j'_1)\cdots(j_t,j'_t)\in\mathcal{R}^j_t$, whenever they are contained in $\mathcal{R}-\{\emptyset_1,\cdots,\emptyset_n\}$.
We will divide our proof into several cases:\vskip5pt

\textbf{Case 1:}\ $l\pm 1\longmapsto l$ or $+(l,l)$.\\
Here we assume that the action of the Kashiwara operator on $\mathbf x$ is given in a way such that $x_q$ satisfies property (a). It means that we either replace $l\pm1$ by $l$ or add the pair $(l,l)$ as described in Section~\ref{section3}. Here we consider again several cases, where each case assumes that a condition described in Theorem~\ref{char2} is violated.\vskip3pt
\textbf{Case 1.1:}\ (1)i) is violated for the pair $(x_{q-1},\theta(x_q))$:\\
Since (1)i) is violated, there must exist an element $\theta(j_{t+\delta_{l,j}-r})=l$ such that \begin{equation}\label{hhj}\sharp\{i_k|1\leq i_k\leq l\}\leq r.\end{equation}
Further, as in Case 1.2 of Theorem~\ref{char}, one can verify that $l$ does not appear in $x_{q-1}$ and if in addition $l\leq i<j$ holds, then $l+1$ occurs in $x_{q-1}$. Consequently $\overline{l}$ must be contained in $\{i'_1,\cdots,i'_s\}$, because otherwise we would obtain a contradiction to the choice of $q$. Hence if we take $i'_p=i'_r=\overline{l}$ it follows immediately that the properties (a), (b) and (c) in (3) hold for the pair $(x_{q-1},x_q)$, which is impossible. For instance (c) is fulfilled with (\ref{hhj}), since $$\sharp\{i_k|1\leq i_k\leq l\}=\sharp\{i_k|i_k< l\}\leq r\leq \sharp\{j_k|j_k<l\}.$$ The other violations of the properties in (1) or (2) can be proven similarly, so that we consider as a next step the following case:
\vskip2pt
\textbf{Case 1.2:}\ (3) is violated for the pair $(\theta(x_q),x_{q+1})$:\\
A simple case-by-case observation shows that this case can only occur if there exists $\theta(i'_p)\geq \theta(i'_r)=\overline{l+1}$, such that (a), (b) and (c) is satisfied. Suppose that there exists an element $\overline{m}$ in the set $\{\theta(i'_{r+1}),\cdots,\theta(i'_p)\}$ $(\overline{m}=\theta(i'_{r+u}))$, such that \begin{equation}\label{mprop}m-\delta_{m,\{i+1,\cdots,n\}}\notin (1-\delta_{m,\{i+1,\cdots,n\}})\{i_1,\ldots,i_s\}\cup \delta_{m,\{i+1,\cdots,n\}}(\mathbf I- \{i'_1,\ldots,i'_s\})\end{equation} and $\overline{m}$ is minimal with this property. In the case where such an element does not exist we set $\overline{m}=\theta(i'_p)$. We claim the following:\vskip 3pt
\textbf{Claim:} Let $\overline{m}$ be as described before, then
\begin{equation}\label{mprop2}u+\min\{l,i\}-\min\{m,i\}-\sharp\{i_k|m\leq i_k<l\}+\sharp\{i'_k|m\leq i'_k<l\}\leq l+1-m.\end{equation}
\textit{Proof of the Claim:}
We consider again various cases, starting with\vskip3pt
$\bullet$ $l\leq i$:\vskip3pt
The minimality of $\overline{m}$ implies $\overline{\theta(i'_{r+1})},\cdots, \overline{\theta(i'_{r+u-1})}\in\{i_1,\cdots,i_s\}$, which yields
$$\underbrace{\sharp\{i'_k|\theta(i'_{r})<i'_k<\overline{m}\}}_{=u-1}\leq \sharp\{i_k|m\leq i_k <l\}
\Longleftrightarrow \sharp\{i'_k|m\leq i'_k <l\}-\sharp\{i_k|m\leq i_k <l\}\leq -u+1.$$\vskip3pt
$\bullet$ $l>i$ and $m\leq i$:\vskip3pt
Using the minimality of $\overline{m}$ we obtain similarly $\overline{\theta(i'_{r+u-1})},\cdots, \overline{\theta(i'_{r+x})}\in\{i_1,\cdots,i_s\}$ and $\overline{\theta(i'_{r+1})}-1,\cdots, \overline{\theta(i'_{r+x-1})}-1\notin\{i'_1,\cdots,i'_s\}$, for some integer $x$. Therefore, the following calculation implies (\ref{mprop2})
$$i-\sharp\{i_k|m\leq i_k\}+\sharp\{i'_k|i'_k <l\}\leq l-(x-1)-\sharp\{i_k|m\leq i_k\}\leq l-\sharp\{i'_k|\theta(i'_{r})<i'_k<\overline{m}\}.$$
$\bullet$ $m>i$:\vskip3pt
The minimality of $\overline{m}$ provides $\overline{\theta(i'_{r+1})}-1,\cdots,\overline{\theta(i'_{r+u-1})}-1\notin\{i'_1,\cdots,i'_s\}$. Hence,
$$\sharp\{i'_k|m\leq i'_k<l\}\leq l-m-\sharp\{i'_k|\theta(i'_{r})<i'_k<\overline{m}\},$$which finishes the proof of the claim.\vskip3pt
As a consequence of (\ref{mprop2}) we obtain 
\begin{flalign*}&\sharp\{j_k|j_k<\overline{\theta(i'_{p})}\}-\sharp\{\theta(i_k)|\theta(i_k)<\overline{\theta(i'_{r})}\}+\sharp\{\theta(i'_k)|\theta(i'_k)<\overline{\theta(i'_{r})}\}-\sharp\{j'_k|j'_k<\overline{\theta(i'_{p})}\}&\\&
\leq\sharp\{j_k|j_k<\overline{i'_{p-\delta_{l,i}}}\}-\sharp\{i_k|m\leq i_k<l\}-\sharp\{i_k| i_k<m\}&\\&+\sharp\{i'_k|m\leq i'_k<l\}+\sharp\{i'_k|i'_k<m\}-\sharp\{j'_k|j'_k<\overline{i'_{p-\delta_{l,i}}}\}
&\\&< l+1-m-\min\{l,i\}+\min\{m,i\}+\max\{0,m-i\}-\max\{0,\overline{i'_{p-\delta_{l,i}}}-j\}-(p-r)
&\\&=l+1-\min\{l,i\}-\max\{0,\overline{i'_{p-\delta_{l,i}}}-j\}-(p-r)
&\\&\leq\max\{0,l+1-i\}-\max\{0,\overline{i'_{p-\delta_{l,i}}}-j\}-(p-r)+\delta_{l,\{1,\cdots,i-1\}},\end{flalign*}
where the first estimation is strict if $l<i$. Hence we have a contradiction to the assumption that (c) holds. \vskip3pt

\textbf{Case 1.3:}\ (3) is violated for the pair $(x_{q-1},\theta(x_q))$:\\
In that case there exists $i'_p=\overline{l+1}\geq i'_r$, such that (a), (b) and (c) is satisfied. Therefore we must have $\overline{l}=i'_{p+1}$, because otherwise we obtain a contradicition to the choice of $q$. It follows
\begin{flalign*}&(p-r)+\sharp\{\theta(j_k)|\theta(j_k)<\overline{i'_{p}}\}-\sharp\{i_k|i_k<\overline{i'_{r}}\}+\sharp\{i'_k|i'_k<\overline{i'_{r}}\}-\sharp\{\theta(j'_k)|\theta(j'_k)<\overline{i'_{p}}\}&\\&
\leq(p+1-r)+\sharp\{j_k|j_k<\overline{i'_{p+1}}\}-\sharp\{i_k|i_k<\overline{i'_{r}}\}+\sharp\{i'_k| i'_k<\overline{i'_{r}}\}-\sharp\{j'_k|j'_k<\overline{i'_{p+1}}\}
&\\&<\max\{0,\overline{i'_r}-i\}-\max\{0,l-j\}\leq \max\{0,\overline{i'_r}-i\}-\max\{0,l+1-j\}+\delta_{l,\{j,\cdots,n\}},\end{flalign*}
where the first estimation is strict whenever $l\geq j$ and thus provides a contradiction to (c).\vskip3pt

\textbf{Case 1.4:}\ (4) is violated for the pair $(\theta(x_q),x_{q+1})$:\\
It is easy to see that this case can never appear.\vskip3pt

\textbf{Case 1.5:}\ (4) is violated for the pair $(x_{q-1},\theta(x_q))$:\\
In that case we have two possibilities, where we start by supposing that there is $i'_p=\overline{l+1}\geq \theta(j'_r)$, such that the properties (a), (b) and (c) are fulfilled. Similar to Case 1.3 we must have an element $\overline{l}=i'_{p+1}$, because otherwise we would obtain a contradicition to the choice of $q$. Subsequently we get
\begin{flalign*}&(i-j)+(t+\delta_{l,j}-s)+(p-r)+\sharp\{\theta(j'_k)|\overline{i'_p}\leq\theta(j'_k)<\overline{\theta(j'_r)}\}-\sharp\{\theta(j_k)|\overline{i'_p}\leq \theta(j_k)<\overline{\theta(j'_r)}\}&\\&
\leq(i-j)+(t+\delta_{l,j}-s)+(p+1-r)+\sharp\{j'_k|l<j'_k\leq\overline{\theta(j'_r)}\}-\sharp\{j_k|l\leq j_k<\overline{\theta(j'_r)}\}
&\\&<\max\{0,\overline{\theta(j'_r)}-j\}-\max\{0,l-j\}\leq \max\{0,\overline{\theta(j'_r)}-j\}-\max\{0,l+1-j\}+\delta_{l,\{j,\cdots,n\}},\end{flalign*}
where the first estimation is strict provided $l\geq j$ meaning that this calculation contradicts once more property (c).\par
The last and second possibility which can occur is: there exists $i'_p\geq \theta(j'_r)=\overline{l+1}$, such that the properties (a), (b) and (c) are fulfilled.\par
Then we make a similar construction as in Case 1.2, namely we suppose that $$\overline{m}\in\{\theta(j'_{r+1}),\ldots,\theta(j'_{t+\delta_{l,j}})\},\mbox{ say }\overline{m}=\theta(j'_{u})=j'_{u-\delta_{l,j}
},$$such that $$m-\delta_{m,\{j+1,\cdots,n\}}\notin (1-\delta_{m,\{j+1,\cdots,n\}})\{j_1,\ldots,j_t\}\ \cup\ \delta_{m,\{j+1,\cdots,n\}}(\mathbf I- \{j'_1,\ldots,j'_t\})$$ and $m\geq \overline{i'_p}$. If such an element exists we choose $m$ maximal with this property and otherwise we set $\overline{m}=i'_p$. Using the maximality, we can verify similar to (\ref{mprop2}) the correctness of
\begin{align*}(u-r)-\sharp\{j_k|m\leq j_k<l\}+\sharp\{j'_k|m\leq j'_k<l\}+\min\{l,j\}-\min\{m,j\}\leq l+1-m.\end{align*}
As a corollary we obtain as usual a contradiction to property (c) (recall that $\theta(x_q)\in\mathcal{R}^j_{t+\delta_{l,j}}$)
\begin{flalign*}&(i-j)+(t+\delta_{l,j}-s)+(p-r)+\sharp\{\theta(j'_k)|\overline{i'_p}\leq\theta(j'_k)<\overline{\theta(j'_r)}\}-\sharp\{\theta(j_k)|\overline{i'_p}\leq \theta(j_k)<\overline{\theta(j'_r)}\}&\\&
\leq(i-j)+(t+\delta_{l,j}-s)+(p-u)+(u-r)+\sharp\{j'_k|\overline{i'_p}\leq j'_k< m\}&\\&+\sharp\{j'_k|m\leq j'_k< l\}-\sharp\{j_k|\overline{i'_p}\leq j_k<m\}-\sharp\{j_k|m\leq j_k<l\}&\\&
<\max\{0,m-j\}-\max\{0,\overline{i'_p}-j\}+l+1-m-\min\{l,j\}+\min\{m,j\}
&\\&\leq l+1-\min\{l,j\}-\max\{0,\overline{i'_p}-j\}\leq \max\{0,l+1-j\}-\max\{0,\overline{i'_p}-j\}+\delta_{l,\{1,\cdots,j-1\}},\end{flalign*}
where the first estimation is strict whenever $l<j$ is satisfied.\vskip5pt

\textbf{Case 2:}\ $\overline{l+1}\longmapsto\overline{l}$.\\
Now we assume that the action of the Kashiwara operator on $\mathbf x$ is given in a way such that $x_q$ satisfies property (b) while (a) is violated, which in particular means that we replace the entry $\overline{l+1}$ by $\overline{l}$. Since the proofs are similar to Case 1, we do not give them in full details. We only consider the case where we presume that property (3) is violated, i.e.\vskip3pt

\textbf{Case 2.1:}\ (3) is violated for the pair $(x_{q-1},\theta(x_q))$:\\
It is easy to see that this case can never appear.\vskip3pt

\textbf{Case 2.2:}\ (3) is violated for the pair $(\theta(x_q),x_{q+1})$:\\
The first possibility which can occur is: there exists $\theta(i'_p)=i'_p\geq \theta(i'_r)=\overline{l}$, such that the properties (a), (b) and (c) are fulfilled. Because of (a) and $i'_r=\overline{l+1}$ we must have $$l+1-\delta_{l+1,\{i+1,\cdots,n\}}\notin (1-\delta_{l+1,\{i+1,\cdots,n\}})\{i_1,\ldots,i_s\}\cup \delta_{l+1,\{i+1,\cdots,n\}}(\mathbf I- \{i'_1,\ldots,i'_s\}),$$ since otherwise we would obtain that $x_q$ satisfies (a) (from Section~\ref{section3}) and thus the Kashiwara operator would act as in Case 1. Accordingly we can apply our assumptions to $i'_p\geq i'_r=\overline{l+1}$ and obtain a contradiction to (c), namely 
\begin{flalign*}&\sharp\{j_k|j_k<\overline{i'_p}\}-\sharp\{\theta(i_k)|\theta(i_k)<l\}+\sharp\{\theta(i'_k)|\theta(i'_k)<l\}-\sharp\{j'_k|j'_k<\overline{i'_p}\}&\\&
=\sharp\{j_k|j_k<\overline{i'_p}\}+\min\{l,i\}-\min\{l,i\}-\sharp\{i_k|i_k<l\}+\sharp\{i'_k|i'_k<l\}-\sharp\{j'_k|j'_k<\overline{i'_p}\}&\\&
<\sharp\{j_k|j_k<\overline{i'_p}\}+\min\{l+1,i\}-\min\{l,i\}-\sharp\{i_k|i_k<l+1\}&\\&+\sharp\{i'_k|i'_k<l+1\}-\sharp\{j'_k|j'_k<\overline{i'_p}\}
&\\&<\min\{l+1,i\}-\min\{l,i\}+\max\{0,l+1-i\}-\max\{0,\overline{i'_p}-j\}-(p-r)&\\&=1+\max\{0,l-i\}-\max\{0,\overline{i'_p}-j\}-(p-r).\end{flalign*}
The second and last possibility which can occur in that case is: there exists $\theta(i'_p)=\overline{l}\geq \theta(i'_r)=i'_r$, such that the properties (a), (b) and (c) are fulfilled. As before we can assume $i'_p=\overline{l+1}$. Suppose that there exists $\overline{m}\in\{i'_{r},\ldots,i'_p\}$ ($\overline{m}=i'_{r+u}$), such that 
$$m-\delta_{m,\{j+1,\cdots,n\}}\notin (1-\delta_{m,\{j+1,\cdots,n\}})\{j_1,\ldots,j_t\}\cup \delta_{m,\{j+1,\cdots,n\}}(\mathbf I- \{j'_1,\ldots,j'_t\})$$ and $m$ is minimal with this property. If such an element does not exist we set $\overline{m}=i'_{r}$. Similar to (\ref{mprop2}) one gets
$$p-(r+u)-\min\{l,j\}+\min\{m,j\}-\sharp\{j_k|l\leq j_k<m\}+\sharp\{j'_k|l\leq j'_k<m\}\leq m-l.$$
Using this inequality we arrive once more at a contradiction, namely
\begin{flalign*}&\sharp\{j_k|j_k<l\}-\sharp\{\theta(i_k)|\theta(i_k)<\overline{i'_r}\}+\sharp\{\theta(i'_k)|\theta(i'_k)<\overline{i'_r}\}-\sharp\{j'_k|j'_k<l\}&\\&
=\sharp\{j_k|j_k<m\}-\sharp\{j_k|l\leq j_k<m\}-\sharp\{i_k|i_k<\overline{i'_r}\}+\sharp\{i'_k|i'_k<\overline{i'_r}\}&\\&-\sharp\{j'_k|j'_k<m\}+\sharp\{j'_k|l\leq j'_k<m\}&\\&
<\max\{0,\overline{i'_r}-i\}-\max\{0,m-j\}-(p-r)-\min\{m,j\}+\min\{l,j\}+m-l&\\&
=\max\{0,\overline{i'_r}-i\}-\max\{0,l-j\}-(p-r).\end{flalign*}
The proof of $\tilde{e}_l\mathbf x\in\mathcal{R}(\lambda)\cup \{0\}$ is similar, which completes the proof.
\endproof
\end{thm}
%
\section{Crystal bases as tuples of integer sequences}\label{section6}
In this section we will verify with Theorem~\ref{mainthm} that the crystal $\mathcal{R}(\lambda)$ can be identified with the crystal graph $B(\lambda)$ obtained from Kashiwara's crystal bases theory. Our strategy here is to show that there exists an isomorphism of $\mathcal{R}(\lambda)$ onto the connected component of $\bigotimes^n_{i=1}B(\omega_i)^{\otimes m_i}$ containing $\bigotimes^n_{i=1}b_i^{\otimes m_i}$, where $b_i$ denotes the highest weight element in $B(\omega_i)$. For the proofs in type $A_n$ we will need a result stated in \cite{K12}, where the affine type A Kirillov-Reshetikhin crystals are realized via polytopes. Especially we will need the realization of level 1 KR-crystals, since they are as classical crystals isomorphic to $B(\omega_i)$. In type $C_n$ we will use a short induction argument to prove our results.

\begin{thm}\label{mainthm}
Let $\lambda$ be an arbitrary dominant integral weight and set $k=\sum_im_i$ as before. Then 
\begin{enumerate} 
\item If $k=1$, we have an isomorphism of crystals $\mathcal{R}(\lambda)\stackrel{\sim}{\longrightarrow}B(\lambda).$
\item If $k>1$ and $j$ is the maximal integer such that $m_j\neq 0$ we obtain a strict crystal morphism $$\phi:\mathcal{R}(\lambda)\longrightarrow \mathcal{R}(\lambda-\omega_j)\otimes \mathcal{R}(\omega_j)$$ mapping $\mathbf x=(x_1,\cdots,x_k,0,0\cdots)$ to the tensor product $\mathbf x_{k-1}\otimes x_k$, where $\mathbf x_{k-1}=(x_1,\cdots,x_{k-1},0,0\cdots).$
\end{enumerate}
\proof 
With the help of Proposition~\ref{prop1} the crystal $\mathcal{R}(\omega_i)$ is characterized as $\bigcup_s\mathcal{R}^i_s$. In the case where $\Lg$ is of type $A_n$ we will consider the map $$x=(i_1,i'_1)\cdots (i_s,i'_s)\longmapsto X,$$ whereby $X$ is a pattern as in \cite[Definition~2.1]{K12}, with filling $$a_{p,q}=\begin{cases} 1,& \text{if $(p,q)=(i_r,i'_r)$ for some $r$}\\ 0,& \text{else.}\end{cases}$$ 
By an inspection of the crystal structure on the KR-crystal $B^{1,i}\cong B(\omega_i)$ (\cite[Section~3.2]{K12}), it is easy to see that this map becomes an isomorphism of crystals.\par
If $\Lg$ is of type $C_n$, the proof of part (1) will proceed by upward induction on $i$. An observation of the crsytal graph of $R(\omega_1)$,
\begin{equation}\label{cga}\emptyset_1\stackrel{1}{\longrightarrow}(1,1)\stackrel{2}{\longrightarrow}\cdots\stackrel{n-1}{\longrightarrow}(1,n-1)\stackrel{n}{\longrightarrow}(1,n)\stackrel{n-1}{\longrightarrow}(1,\overline{n-1})\stackrel{n-2}{\longrightarrow}\cdots\stackrel{1}{\longrightarrow}(1,\overline{1})\end{equation}
proves the initial step. Now we assume the correctness of the claim for all integers less than $i$, especially we have $\mathcal{R}(\omega_{i-1})\cong B(\omega_{i-1})$. For the purpose of completing the induction we consider the injective map $$\eta:\mathcal{R}(\omega_{i})\longrightarrow \mathcal{R}(\omega_{i-1})\otimes \mathcal{R}(\omega_{1})$$ given by 
$$\eta(x)=\begin{cases}\emptyset_{i-1}\otimes(1,i-1),& \text{if s=0}\\
\emptyset_{i-1}\otimes(1,i'_1),& \text{if $s=1,i_1=i$}\\ (i_1,i-1)^{(1-\delta_{i_1,i})}(i_2,i'_1)\cdots(i_s,i'_{s-1})\otimes(1,i'_s),& \text{else.}\end{cases}$$ 
If $\eta$ would be a strict crystal morphism, we would get $\mathcal{R}(\omega_{i})\cong \mathcal{I}m(\eta)\cong B(\omega_i)$, which finishes the induction. Therefore we prove that $\eta$ is a strict crystal morphism, where we consider the cases $s=0$ and $s=1,i_1=i$ separately. In the separated cases we draw a part of the crystal graph in order to see that the properties of a crystal morphism hold.
If $s=0$ we obtain
$$\emptyset_i\stackrel{i}{\longrightarrow}(i,i),\quad \emptyset_{i-1}\otimes (1,i-1)=\eta(\emptyset_i)\stackrel{i}{\longrightarrow}\emptyset_{i-1}\otimes(1,i)=\eta((i,i))$$
and if $s=1,i_1=i$ we get

\begin{enumerate}
  \item if $i'_1\leq n-1$
  $$\begin{xy}
  \xymatrix@C=0pt{
  &\emptyset_i^{\delta_{i'_1,i}}(i,i'_1-1)^{(1-\delta_{i'_1,i})}\ar[d]^{i'_1}\\
     &(i,i'_1) \ar[dl]^{i'_1+1} \ar[rd]_{i-1}& \\
                    (i,i'_1+1) & &(i-1,i'_1) 
  }
\end{xy}\quad \begin{xy}
  \xymatrix@C=0pt{
  &\emptyset_{i-1}\otimes(1,i'_1-1)\ar[d]^{i'_1}\\
     &\emptyset_{i-1}\otimes(1,i'_1) \ar[dl]^{i'_1+1} \ar[rd]_{i-1}& \\
                    \emptyset_{i-1}\otimes(1,i'_1+1) & &(i-1,i-1)\otimes(1,i'_1) 
  }
\end{xy}$$
\item if $i'_1\geq n$ and $i'_1\neq \overline{i}$
$$\begin{xy}
  \xymatrix{
  &(i,\overline{\overline{i'_1}+1})\ar[d]^{\overline{i'_1}}\\
     &(i,i'_1) \ar[dl]^{\overline{i'_1}-1} \ar[rd]_{i-1}& \\
                    (i,\overline{\overline{i'_1}-1}) & &(i-1,i'_1) 
  }
\end{xy}\quad \begin{xy}
  \xymatrix@C=0pt{
  &\emptyset_{i-1}\otimes(1,\overline{\overline{i'_1}+1})\ar[d]^{\overline{i'_1}}\\
     &\emptyset_{i-1}\otimes(1,i'_1) \ar[dl]^{\overline{i'_1}-1} \ar[rd]_{i-1}& \\
                    \emptyset_{i-1}\otimes(1,\overline{\overline{i'_1}-1}) & &(i-1,i-1)\otimes(1,i'_1) 
  }
\end{xy}$$
  \item if $i'_1\geq n$ and $i'_1= \overline{i}$
$$\begin{xy}
  \xymatrix{
  &(i,\overline{i+1})\ar[d]^{i}\\
     &(i,\overline{i})\ar[d]_{i-1}& \\
                  &(i-1, \overline{i})\ar[d]_{i-1}&  \\
                  &(i-1, \overline{i-1})
  }
\end{xy}\quad \begin{xy}
  \xymatrix{
  &\emptyset_{i-1}\otimes(1,\overline{i+1})\ar[d]^{i}\\
     &\emptyset_{i-1}\otimes(1,\overline{i})\ar[d]_{i-1}& \\
                  &(i-1,i-1)\otimes(1,\overline{i}) \ar[d]_{i-1}&  \\
                  &(i-1,i-1)\otimes(1, \overline{i-1})
  }
\end{xy}$$
 \end{enumerate}

Thus from now on we can presume that $s>1$ or $s=1,i_1\neq i$. Let $l\in I$ be an arbitrary integer. We set for convenience $\eta(x)=x_1\otimes x_2$, then \begin{align*}\varphi_l(\eta(x))=\max\{\varphi_l(x_2),\varphi_l(x_1)+\varphi_l(x_2)-\epsilon_l(x_2)\}=\begin{cases}
\varphi_l(x_1)+\varphi_l(x_2),& \text{if $l\neq i'_s,\overline{i'_s}$}\\
\max\{0,\varphi_l(x_1)-1\},& \text{else.}\end{cases}\end{align*}

The proof of $\varphi_l(\eta(x))=\varphi_l(x)$ is an intensive investigation of the properties (a)-(d) listed in Section~\ref{section3}. To avoid confusion with indices we consider only the following case:\vskip 3pt

$\bullet$ $l=i'_s$ or $l=\overline{i'_s}$:\vskip 3pt
The case $l=i'_s$ is very simple, because $x$ satisfies neither (a) nor (b), which yields $\varphi_l(x)=0$. Furthermore, since $\overline{l+1}=\overline{i'_s+1}\notin\{i'_1,\cdots,i'_s\}$ we get $\varphi_l(x_1)\leq 1$. Consequently $\varphi_l(\eta(x))=0$. So it remains to consider the case $l=\overline{i'_s}$.
Here we claim the following:\par
let $(\tilde{a})$ be the property which arises from $(a)$ by erasing the condition $\overline{l}\notin\{i'_1,\cdots,i'_s\}$, then $x$ satisfies $(\tilde{a})$ if and only if $x_1$ satisfies $(\tilde{a})$.\par
First we observe that the property $\overline{l+1}\notin\{i'_1,\cdots,i'_s\}$ is equivalent to $\overline{l+1}\notin\{i'_1,\cdots,i'_{s-1}\}$, so that we can ignore it.
We start the proof by assuming that $x_1$ satisfies $(\tilde{a})$ and $x$ does not satisfy $(\tilde{a})$, e.g. $l\in\{i_1,\cdots,i_s,i'_1,\cdots,i'_{s}\}$. This case is only possible if $l=i_1=i$ and thus $l-1=i-1\in\{i'_1,\cdots,i'_{s-1}\}$, which is a contradiction. Since $l-1\in\{i'_1,\cdots,i'_{s}\}$ is always fulfilled if $l>i>i-1$ we can assume that $l<i$ and additionaly $l+1\notin\{i_1,\cdots,i_s\}.$ Without loss of generality we set $l=i-1$, because else we have $l+1\in\{i_1^{(1-\delta_{i_1,i})},i_2,\cdots,i_s\}$. As a consequence we get $$l+1=i\notin\{i_1,\cdots,i_s\}\Longrightarrow i_1\leq i-1\Longrightarrow l\neq i-1,$$ which is again a contradiction.\par 
According to these calculations, $x$ must satisfy $(\tilde{a})$. To show the other direction let $x_1$ violate one of the properties in $(\tilde{a})$, e.g. $l\in\{(i_1,i-1)^{(1-\delta_{i_1,i})},i_2,\cdots,i_s,i'_1,\cdots,i'_{s-1}\}$. Then we must have $l=i-1$ and $i_1<i$ which violates $l+1\in\{i_1,\cdots,i_s\}.$ The only additional possibility which can occur, such that $x_1$ does not satisfy $(\tilde{a})$, is $l<i-1$ and $l+1\notin\{i_1^{(1-\delta_{i_1,i})},i_2,\cdots,i_s\}$. Then we get $l+1=i_1=i$, which is a contradiction to $l<i-1$ and so the reverse direction is also completed.\par
According to this we have $\varphi_l(x)=1\Longrightarrow \varphi_l(x_1)=2$ and $\varphi_l(x)=0\Longrightarrow \varphi_l(x_1)\leq 1$.\vskip 3pt

To show the existence of a morphism of crystals we have to show (among others) the weight invariance of $\eta$, which is proven by the following calculation:
\begin{align}\label{winv}\wt(\eta(x))&=\omega_{i-1}+\omega_1-\sum^{s-1}_{j=1}\alpha_{i_{j+1},i'_j}-\alpha_{1,i'_s}-(1-\delta_{i_1,i})\alpha_{i_1,i-1}&\\&\notag=\omega_{i-1}+\omega_1-\alpha_{1,i-1}-\sum^{s}_{j=1}\alpha_{i_{j},i'_j}=\omega_i-\sum^{s}_{j=1}\alpha_{i_{j},i'_j}.\end{align}
The verification of $\epsilon_l(\eta(x))=\epsilon_l(x)$ is therefore proven with Definition~\ref{abcrystal} (1) and (\ref{winv}).
Suppose now that $\tilde{f_l}\eta(x)\neq 0$ and $\tilde{f}_l$ acts on the second tensor. A short investigation of the crystal graph (\ref{cga}) yields \begin{equation}\label{ggur}l=\begin{cases}i'_s+1,& \text{if $i'_s\leq n-1$}\\ 
\overline{i'_s}-1,& \text{else.}\end{cases}\end{equation} Since $\varphi_l(\eta(x))=\varphi_l(x)=1\neq 0$, we obtain that $$\tilde{f}_lx=\begin{cases}(i_1,i'_1)\cdots(i_s,i'_s+1),& \text{if $l=i'_s+1$}\\ 
(i_1,i'_1)\cdots(i_s,\overline{\overline{i'_s}-1}),& \text{$l=\overline{i'_s}-1$,}\end{cases}$$ because if $l=\overline{i'_s}-1$ and $x$ would satisfy (a), it would automatically follow that $x$ satisfies (a') which is impossible since $\varphi_l(x)=1\neq 2$. Thus $\tilde{f_l}\eta(x)=\eta(\tilde{f_l}x).$\par
If $\tilde{f_l}\eta(x)\neq 0$ and $\tilde{f}_l$ acts on the first tensor we get with the tensor product rule $\varphi_l(x_1)>\epsilon_l(x_2)$. Note that the operation with the Kashiwara operator $\tilde{f}_l$ would change the entry $i'_s$ in $x$ if and only if either $l=i'_s+1$ or $l=\overline{i'_s}-1$ and $x$ is not subject to (a). Our goal is to show here that $\tilde{f}_l$ has no effect on $i'_s$. If $l=i'_s+1$ we would get $\varphi_l(x_1)=0$ and thus a contradicition to $\varphi_l(x_1)>\epsilon_l(x_2)$. In the case where $l=\overline{i'_s}-1$, we have that (b) is not fulfilled for $x_1$. Therefore $x_1$ must be subject to (a). 
Consequently we obtain that $x$ must fulfill (a) as well, because otherwise we would end in a contradiction, namely:\par
the only property in (a) which can be violated is $l\notin\{i_1,\cdots,i_s,i'_1,\cdots,i'_s\}.$ So if $l$ is contained in the aforementioned set we get by the definition of $\eta$ that $l=i_1=i$. Hence $l-1=i-1$ must be contained in the set $\{i\leq i'_1,\cdots,i'_{s-1}\}$, which is impossible.\par
Thus, the entry $i'_s$ stays unchanged in $\tilde{f}_lx$ which provides $\tilde{f_l}\eta(x)=\eta(\tilde{f_l}x).$
The proof of $\tilde{e_l}\eta(x)=\eta(\tilde{e_l}x)$ is similar, which completes the proof of (1).\vskip 3pt
In order to prove (2) we will check as in (1) step for step the properties of a morphism of crystals.
We get \begin{align*}\varphi_l(\mathbf x_{k-1}\otimes x_k)&=\max\{\varphi_l(x_k),\varphi_l(\mathbf x_{k-1})+\varphi_l(x_k)-\epsilon_l(x_k)\}&\\&=\max\{\varphi_l(x_k),\sigma^{k}_l(\mathbf x_{k-1})-\tau^{k-1}_l(\mathbf x_{k-1})&\\&
-\min\{\sigma^j_l(\mathbf x_{k-1})-\tau^j_l(\mathbf x_{k-1})|1\leq j\leq k-1\}+\varphi_l(x_k)-\epsilon_l(x_k)\}&\\&
=\max\{\varphi_l(x_k),\varphi_l(x_k)+\sigma^k_l(\mathbf x)-\tau^k_l(\mathbf x)-\min\{\sigma^j_l(\mathbf x)-\tau^j_l(\mathbf x)|1\leq j\leq k-1\}\}&\\&
=\sigma^{k+1}_l(\mathbf x)-\tau^k_l(\mathbf x)-\min\{\sigma^j_l(\mathbf x)-\tau^j_l(\mathbf x)|1\leq j\leq k\}&\\&=\varphi_l(\mathbf x).\end{align*}
The same is trivially fulfilled for $\epsilon_l$ because of Definition~\ref{abcrystal} (1) and the weight invariance of $\phi$. Now suppose that $f_l(\mathbf x)=p$ is as in (\ref{1}). If we apply the Kashiwara operator $\tilde{f}_l$ to the tensor product $\mathbf x_{k-1}\otimes x_k$ we obtain with the above calculations \begin{align*}\tilde{f}_{l}(\mathbf x_{k-1}\otimes x_k)&=\begin{cases} \tilde{f}_{l}\mathbf x_{k-1}\otimes x_k, \text{ if $\varphi_l(\mathbf x_{k-1})>\epsilon_l(x_k)$}\\
\mathbf x_{k-1}\otimes \tilde{f}_{l}x_k, \text{ if $\varphi_l(\mathbf x_{k-1})\leq \epsilon_l(x_k)$}\end{cases}&\\&=\begin{cases} \tilde{f}_{l}\mathbf x_{k-1}\otimes x_k, \text{ if $\sigma^k_l(\mathbf x)-\tau^k_l(\mathbf x)>\sigma^p_l(\mathbf x)-\tau^p_l(\mathbf x)$}\\
\mathbf x_{k-1}\otimes \tilde{f}_{l}x_k, \text{ if $\sigma^k_l(\mathbf x)-\tau^k_l(\mathbf x)\leq \sigma^p_l(\mathbf x)-\tau^p_l(\mathbf x)$}.\end{cases}\end{align*}
According to this we get that $p\in\{1,\ldots,k-1\}$ if $\tilde{f}_l$ acts on the first tensor and $p=k$ if $\tilde{f}_l$ acts on the second tensor. The proof for the Kashiwara operator $\tilde{e}_l$ works similar.
\endproof
\end{thm}
\begin{cor}\label{mainthm1}
We have an isomorphism of crystals $$\mathcal{R}(\lambda)\cong B(\lambda).$$
\proof
The proof will proceed by induction on $k$, where the initial step is exactly part (1) of Theorem~\ref{mainthm}. If $k>1$ and $j$ is the maximal integer where $m_j$ is non zero, we can assume with the induction hypothesis that $\mathcal{R}(\lambda-\omega_j)\cong B(\lambda-\omega_j)$ and $\mathcal{R}(\omega_j)\cong B(\omega_j)$. The rest of the proof is done with part (2) of Theorem~\ref{mainthm}, since the map $\phi$ is injective and the image is a connected crystal containing the highest weight element $r_{\lambda-\omega_j}\otimes r_{\omega_j}.$
\endproof
\end{cor}
%

\bibliographystyle{plain}
\bibliography{integersequences-crystal-biblist}
\end{document}